\documentclass{article}
\usepackage{arxiv}
\usepackage[utf8]{inputenc}
\usepackage[T1]{fontenc} 

\usepackage{amsmath,amssymb,amsfonts,amsthm, graphicx, lscape, booktabs, pdflscape, subcaption, longtable, float, enumitem, algorithm, algpseudocode, hyperref, siunitx, url, float, microtype, graphicx, natbib, doi, rotating}

\graphicspath{{figures/}{runs/ccd_tf_validation_full/figures/}}
\numberwithin{equation}{section}
\theoremstyle{definition}
\newtheorem{assumption}{Assumption}
\newtheorem{definition}{Definition}
\newtheorem{proposition}{Proposition}
\newtheorem{remark}{Remark}

\newtheorem{corollary}{Corollary}
\newtheorem{theorem}{Theorem}
\newtheorem{example}{Example}

\title{Core-Conditioned Regularized Matrix Tri-Factorization for High-Dimensional Structured Systems}

\author{Ronald Katende 
Department of Mathematics\\
Kabale University\\
Kikungiri Hill, Katuna Road, 317, Kabale, Uganda\\
\texttt{rkatende92@gmail.com} \\
\\
}

\date{}

\renewcommand{\headeright}{}
\renewcommand{\undertitle}{}
\renewcommand{\shorttitle}{Core-Conditioned Regularized Matrix Tri-Factorization for High-Dimensional Structured Systems}

\hypersetup{
pdftitle={Core-Conditioned Regularized Matrix Tri-Factorization for High-Dimensional Structured Systems},
pdfauthor={Ronald Katende},
pdfkeywords={regularized matrix factorization, low-rank approximation, matrix tri-factorization, core conditioning, structured matrices, numerical stability, nonconvex optimization},
}

\begin{document}
\maketitle

\begin{abstract}
	This paper studies a regularized matrix tri-factorization \(A\approx PDQ\), where \(P\) and \(Q\) are side factors and \(D\) is a central core whose conditioning can be explicitly regularized or constrained. The formulation is a structured low-rank approximation framework, not a replacement for LU, QR, Cholesky, or the singular value decomposition. In the unregularized full-data Frobenius rank-\(r\) problem, truncated SVD remains the optimal benchmark. The contribution here concerns the regularized and core-conditioned setting, where reconstruction accuracy is treated together with factor scale, numerical conditioning, perturbation behavior, and weighted approximation.
	
	The analysis establishes the algebraic scope of the \(PDQ\) representation, proves existence of minimizers under coercive regularization, identifies the non-uniqueness induced by latent-space transformations, derives well-posed block updates for the quadratic full-data objective, and gives product-level perturbation bounds. For exact alternating minimization in the full-data quadratic case, it proves descent, boundedness of iterates, and convergence to a critical point under standard Kurdyka--\L{}ojasiewicz assumptions. A full multi-seed validation indicates competitive behavior in noisy and ill-conditioned low-rank approximation while reporting diagnostics not provided by purely spectral baselines, including the learned core condition number and block-system conditioning. The validation also clarifies the method's limits: randomized SVD remains faster for pure spectral compression, and the current weighted missing-entry variant is not uniformly competitive with matrix-completion baselines. The framework is therefore best viewed as a regularized and diagnostically transparent tri-factorization for settings where approximation quality and numerical conditioning must be controlled jointly.
\end{abstract}

\keywords{regularized matrix factorization; low-rank approximation; matrix tri-factorization; core conditioning; structured matrices; numerical stability; nonconvex optimization.}

\section{Introduction}

High-dimensional matrices are central objects in numerical linear algebra, scientific computing, optimization, statistics, signal processing, network analysis, inverse problems, and machine learning. They may arise as discretized differential operators, covariance matrices, kernel matrices, graph adjacency matrices, recommender-system arrays, image and signal matrices, or partially observed data tables. In these settings, the computational goal is often not merely to factorize a matrix exactly. One often seeks a representation that is lower-dimensional, numerically stable, interpretable, and compatible with structural information in the data or operator.

Classical matrix factorizations provide the foundation for this subject. LU, QR, Cholesky, and singular value decompositions are standard tools for solving linear systems, least-squares problems, eigenvalue problems, and low-rank approximation tasks \citep{GolubVanLoan2013}. Their numerical stability and implementation issues are treated systematically in modern numerical linear algebra \citep{Higham2002}. The singular value decomposition is especially important because it gives an optimal rank-\(r\) approximation in the classical unweighted Frobenius-norm setting \citep{EckartYoung1936}. Mirsky later extended this optimality viewpoint to unitarily invariant norms \citep{Mirsky1960}.

Despite their foundational role, classical decompositions do not directly solve every modern high-dimensional approximation problem. Direct factorizations may be expensive for dense large-scale systems, while purely algebraic decompositions usually do not encode sparsity, missingness, factor regularity, conditioning constraints, or application-dependent structure inside the factorization objective. These issues are especially visible in data matrices, where approximate low-rank structure often coexists with noise, missing entries, anisotropic scaling, ill-conditioning, and heterogeneous sampling patterns \citep{UdellTownsend2019}.

Several research directions address these challenges. Randomized numerical linear algebra reduces the cost of matrix approximation by sketching, sampling, and randomized projection \citep{HalkoMartinssonTropp2011}. A later survey by Martinsson and Tropp gives a broad account of randomized methods for low-rank approximation, regression, kernel approximation, and linear systems \citep{MartinssonTropp2020}. CUR decompositions improve interpretability by approximating a matrix using selected columns and rows \citep{MahoneyDrineas2009}. Deterministic CUR constructions also provide scalable alternatives to purely spectral approximations \citep{BoutsidisMahoneyDrineas2008}. Rank-revealing QR methods form another classical route to structure-aware low-rank approximation \citep{Chan1987}.

Regularized and structured factorization methods form a second major line of work. Nonnegative matrix factorization imposes positivity and has been widely used for parts-based representations \citep{LeeSeung1999}. Sparse principal component analysis introduces sparsity into principal directions \citep{ZouHastieTibshirani2006}. Penalized matrix decomposition provides a flexible framework for sparse PCA and canonical correlation analysis \citep{WittenTibshiraniHastie2009}. Dictionary learning and sparse coding similarly formulate representation learning as a regularized matrix factorization problem \citep{MairalBachPonceSapiro2010}. Generalized low-rank models extend matrix factorization beyond squared loss by allowing different losses and regularizers adapted to the data type \citep{UdellHornZadehBoyd2016}.

Convex low-rank methods provide another important reference point. Nuclear-norm minimization gives a tractable convex surrogate for rank minimization \citep{Fazel2002}. Exact matrix completion via convex optimization established a central recovery result for low-rank matrices observed through a subset of entries \citep{CandesRecht2009}. The noisy matrix completion setting was analyzed by Cand\`es and Plan \citep{CandesPlan2010}. Recht, Fazel, and Parrilo gave a broad treatment of minimum-rank recovery through nuclear-norm minimization \citep{RechtFazelParrilo2010}.

Nonconvex factorized formulations are also now standard in scalable low-rank optimization. The Burer--Monteiro approach represents positive semidefinite variables by low-rank factors and solves the resulting nonconvex problem directly \citep{BurerMonteiro2003}. Optimization algorithms on matrix manifolds provide a geometric framework for low-rank constrained problems \citep{AbsilMahonySepulchre2008}. Alternating minimization has been analyzed for low-rank matrix completion under suitable incoherence and sampling assumptions \citep{JainNetrapalliSanghavi2013}. More recent work has clarified the relation between factorized and manifold formulations in low-rank matrix optimization \citep{LuoLiZhang2023}. Nuclear-norm regularized low-rank optimization can also be accelerated through Burer--Monteiro decompositions \citep{LeeLiangTangToh2024}.

This paper belongs to the optimization-based matrix factorization tradition. We study a regularized tri-factorization of the form
\[
A \approx P D Q,
\qquad
P\in\mathbb{R}^{m\times r},\quad
D\in\mathbb{R}^{r\times r},\quad
Q\in\mathbb{R}^{r\times n},
\]
where \(r\leq \min\{m,n\}\). The factor \(D\) is a low-dimensional core, while \(P\) and \(Q\) are left and right embedding factors. Since every rank-\(r\) matrix can be represented in such a form, the novelty is not the algebraic identity \(PDQ\) by itself. The contribution is the regularized and core-conditioned use of this representation, where the factors are obtained by minimizing a structured approximation objective.

The basic model is
\[
\min_{P,D,Q}
\frac12\|A-PDQ\|_F^2
+
\mathcal{R}(P,D,Q),
\]
where \(\mathcal{R}\) controls factor scale, sparsity, smoothness, normalization, conditioning, or other structural properties. For partially observed data, the residual can be replaced by
\[
\frac12\|W\odot(A-PDQ)\|_F^2,
\]
where \(W\) is a nonnegative weight or observation mask and \(\odot\) denotes the Hadamard product. This places the method within regularized low-rank approximation, but with explicit attention to the role of the central core \(D\).

The distinction from truncated SVD is essential. For a fully observed matrix and the unconstrained rank-\(r\) Frobenius approximation problem, truncated SVD is optimal \citep{EckartYoung1936}. The present method is not proposed as a replacement for SVD in that classical setting. It is intended for regimes where the approximation problem includes additional requirements such as conditioning control, weighted or missing observations, structural regularization, numerical robustness, or downstream stability constraints.

The main contributions are as follows.

\begin{enumerate}
\item We formulate a core-conditioned regularized tri-factorization framework for structured low-rank matrix approximation. The formulation accommodates full-data and weighted-data objectives, but the strongest theoretical results in this paper are stated for the full-data quadratic-regularized case.

\item We state the algebraic scope of the factorization precisely. The model represents matrices of rank at most \(r\), and therefore should be understood as a regularized parametrization of the rank-constrained matrix class rather than a new algebraic rank class.

\item We prove existence of minimizers under explicit coercivity assumptions. This avoids the incorrect compactness argument that would arise if the unconstrained factor spaces were treated as compact.

\item We identify the non-uniqueness induced by latent-space transformations and separate product-level approximation from factor-level identifiability. This prevents unsupported uniqueness claims for the factors \(P,D,Q\).

\item We derive closed-form block updates for the quadratic-regularized full-data objective and give dense and sparse per-sweep cost bounds. These bounds are per-iteration bounds, not total-runtime guarantees.

\item We analyze the role of the core matrix \(D\) by giving explicit core-conditioning and block-system bounds. The conditioning of \(D\) is controlled through the admissible core model or regularization, not inherited automatically from \(A\).
\end{enumerate}

The representation \(PDQ\) is related to matrix tri-factorization, regularized low-rank factorization, and nonconvex factorized optimization. The distinction in this paper is not the existence of a three-factor representation. The distinction is the explicit treatment of the central factor \(D\) as a core whose conditioning can be monitored, regularized, or constrained, together with product-level perturbation bounds and block-system conditioning estimates.

The proposed method should therefore be viewed as a regularized, core-conditioned low-rank approximation tool. It is not a universal substitute for LU, QR, Cholesky, SVD, Krylov methods, randomized decompositions, or matrix-completion algorithms. Each of these methods remains appropriate in its own setting. The contribution of this work is to provide an optimization-based tri-factorization in situations where low-rank approximation, factor regularization, and numerical conditioning must be handled jointly.

\section{Preliminaries and Problem Formulation}

Let \(A\in\mathbb{R}^{m\times n}\). We denote the Frobenius norm by \(\|\cdot\|_F\), the spectral norm by \(\|\cdot\|_2\), and the singular values of \(A\) by
\[
\sigma_1(A)\geq \sigma_2(A)\geq \cdots \geq 0.
\]
For a nonsingular square matrix \(B\), its spectral condition number is
\[
\kappa(B)=\frac{\sigma_{\max}(B)}{\sigma_{\min}(B)}.
\]
If \(B\) is singular, we set \(\kappa(B)=+\infty\). This convention is standard in numerical stability analysis \citep{Higham2002}.

Let \(r\leq \min\{m,n\}\). Define the tri-factorized model class
\[
\mathcal{M}_r
=
\left\{
PDQ:
P\in\mathbb{R}^{m\times r},
D\in\mathbb{R}^{r\times r},
Q\in\mathbb{R}^{r\times n}
\right\}.
\]
Every element of \(\mathcal{M}_r\) has rank at most \(r\). Conversely, every matrix of rank at most \(r\) belongs to \(\mathcal{M}_r\). Hence the model class is exactly the usual rank-constrained matrix class, written in a form that allows separate regularization of the left factor, core factor, and right factor.

\begin{definition}[Core-conditioned regularized tri-factorization]
Let \(A\in\mathbb{R}^{m\times n}\), and let \(r\leq \min\{m,n\}\). A core-conditioned regularized tri-factorization of \(A\) is any solution of
\begin{equation}
\min_{P,D,Q}
F(P,D,Q)
:=
\frac12\|A-PDQ\|_F^2
+
\mathcal{R}(P,D,Q),
\label{eq:main_objective}
\end{equation}
where
\[
P\in\mathbb{R}^{m\times r},\qquad
D\in\mathbb{R}^{r\times r},\qquad
Q\in\mathbb{R}^{r\times n},
\]
and \(\mathcal{R}\) is a proper lower semicontinuous regularizer.
\end{definition}

A basic quadratic regularizer is
\begin{equation}
\mathcal{R}(P,D,Q)
=
\frac{\alpha_P}{2}\|P\|_F^2
+
\frac{\alpha_D}{2}\|D\|_F^2
+
\frac{\alpha_Q}{2}\|Q\|_F^2,
\qquad
\alpha_P,\alpha_D,\alpha_Q>0.
\label{eq:quadratic_regularizer}
\end{equation}
Additional structure may be imposed through sparsity penalties, nonnegativity constraints, normalization constraints, or a core-conditioning penalty. Sparse and penalized factor models are standard in statistical learning and signal processing \citep{ZouHastieTibshirani2006}. Penalized matrix decomposition is one widely used instance of this principle \citep{WittenTibshiraniHastie2009}.

To control the numerical stability of the core, one may impose an explicit feasible set
\begin{equation}
\mathcal{D}_{\delta,M}
=
\left\{
D\in\mathbb{R}^{r\times r}:
\sigma_{\min}(D)\geq \delta,\ 
\|D\|_2\leq M
\right\},
\qquad
0<\delta\leq M<\infty.
\label{eq:core_constraint}
\end{equation}
Then every feasible \(D\) satisfies
\[
\kappa(D)\leq \frac{M}{\delta}.
\]
This makes the conditioning assumption explicit rather than hidden.

For missing or weighted observations, we use the objective
\begin{equation}
F_W(P,D,Q)
=
\frac12\|W\odot(A-PDQ)\|_F^2
+
\mathcal{R}(P,D,Q),
\label{eq:weighted_objective}
\end{equation}
where \(W\in\mathbb{R}^{m\times n}\) is a nonnegative weight matrix. Matrix completion and weighted low-rank approximation motivate this type of formulation \citep{CandesRecht2009}. Noisy matrix completion motivates the same formulation when the observed entries are contaminated by measurement error \citep{CandesPlan2010}.

\subsection{Algebraic Scope}

The following proposition states exactly what the factorization can and cannot represent.

\begin{proposition}[Exact representability at sufficient rank]
Let \(A\in\mathbb{R}^{m\times n}\) have rank \(s\). If \(r\geq s\), then there exist matrices
\[
P\in\mathbb{R}^{m\times r},\qquad
D\in\mathbb{R}^{r\times r},\qquad
Q\in\mathbb{R}^{r\times n}
\]
such that \(A=PDQ\). If \(r<s\), then exact representation is impossible in general.
\end{proposition}

\begin{proof}
Let \(A=U_s\Sigma_sV_s^\top\) be a compact singular value decomposition of \(A\), where \(U_s\in\mathbb{R}^{m\times s}\), \(\Sigma_s\in\mathbb{R}^{s\times s}\), and \(V_s\in\mathbb{R}^{n\times s}\). If \(r=s\), choose
\[
P=U_s,\qquad D=\Sigma_s,\qquad Q=V_s^\top.
\]
Then \(A=PDQ\). If \(r>s\), the same construction is embedded into \(r\) dimensions by adding zero columns and rows. Conversely, for any conformable matrices \(P,D,Q\), the product \(PDQ\) has rank at most \(r\). Therefore a matrix with rank greater than \(r\) cannot be represented exactly by such a product.
\end{proof}

This proposition is important for the interpretation of the method. The factorization does not enlarge the rank-\(r\) model class. Its purpose is to regularize and condition the search within that class.

\subsection{Existence of Minimizers}

The factor spaces in \eqref{eq:main_objective} are not compact. Therefore existence of a minimizer does not follow from compactness unless explicit compact constraints are imposed. A correct existence result follows from coercivity and lower semicontinuity, which are standard assumptions in finite-dimensional variational analysis \citep{RockafellarWets1998}.

\begin{assumption}[Coercive regularization]
The regularizer \(\mathcal{R}\) is proper, lower semicontinuous, bounded from below, and coercive. That is,
\[
\mathcal{R}(P,D,Q)\to+\infty
\]
whenever
\[
\|P\|_F+\|D\|_F+\|Q\|_F\to+\infty.
\]
\end{assumption}

\begin{proposition}[Existence of a global minimizer]
Let \(A\in\mathbb{R}^{m\times n}\). Suppose \(\mathcal{R}\) satisfies the coercive regularization assumption. Then the objective
\[
F(P,D,Q)
=
\frac12\|A-PDQ\|_F^2
+
\mathcal{R}(P,D,Q)
\]
has at least one global minimizer.
\end{proposition}

\begin{proof}
The residual term is continuous and nonnegative. Since \(\mathcal{R}\) is bounded from below, the objective is bounded from below. Since \(\mathcal{R}\) is coercive, every minimizing sequence for \(F\) is bounded. In finite-dimensional Euclidean space, every bounded sequence has a convergent subsequence. The residual term is continuous and \(\mathcal{R}\) is lower semicontinuous, so \(F\) is lower semicontinuous. Therefore the limit of a convergent minimizing subsequence attains the infimum.
\end{proof}

The quadratic regularizer in \eqref{eq:quadratic_regularizer} is coercive when \(\alpha_P,\alpha_D,\alpha_Q>0\). Hence the existence result applies directly to the basic model.

\subsection{Non-Uniqueness and Identifiability}

The individual factors \(P,D,Q\) are not identifiable without additional constraints. For any invertible matrices \(S,T\in\mathbb{R}^{r\times r}\),
\begin{equation}
PDQ
=
(PS)(S^{-1}DT^{-1})(TQ).
\label{eq:gauge_transform}
\end{equation}
Thus, infinitely many triples may represent the same product. This type of non-uniqueness is common in factorized low-rank optimization \citep{BurerMonteiro2003}. It is also one reason why manifold formulations and factorized formulations require careful comparison \citep{LuoLiZhang2023}.

The identifiable object in the unconstrained model is the product
\[
X=PDQ,
\]
not necessarily the individual factors. Factor-level identifiability requires additional gauge fixing, such as diagonal \(D\), column normalization, sign conventions, ordering constraints, sparsity assumptions, or orthogonality constraints.

\subsection{Quadratic Block Updates}

Consider the quadratic objective
\begin{equation}
F(P,D,Q)
=
\frac12\|A-PDQ\|_F^2
+
\frac{\alpha_P}{2}\|P\|_F^2
+
\frac{\alpha_D}{2}\|D\|_F^2
+
\frac{\alpha_Q}{2}\|Q\|_F^2,
\label{eq:quadratic_objective}
\end{equation}
where \(\alpha_P,\alpha_D,\alpha_Q>0\). Holding two blocks fixed, each update is a regularized least-squares problem.

For fixed \(D\) and \(Q\), the update for \(P\) is
\begin{equation}
P
=
A Q^\top D^\top
\left(
DQQ^\top D^\top+\alpha_P I_r
\right)^{-1}.
\label{eq:P_update}
\end{equation}
For fixed \(P\) and \(D\), the update for \(Q\) solves
\begin{equation}
\left(
D^\top P^\top P D+\alpha_Q I_r
\right)Q
=
D^\top P^\top A.
\label{eq:Q_update}
\end{equation}
For fixed \(P\) and \(Q\), the update for \(D\) solves
\begin{equation}
(P^\top P)D(QQ^\top)+\alpha_DD
=
P^\top A Q^\top.
\label{eq:D_update}
\end{equation}
The \(P\)- and \(Q\)-updates are well posed because the coefficient matrices are positive definite when the corresponding regularization parameters are positive. The \(D\)-update is also uniquely solvable because the associated linear operator has eigenvalues
\[
\lambda_i(P^\top P)\lambda_j(QQ^\top)+\alpha_D,
\]
which are strictly positive for \(\alpha_D>0\).

The alternating scheme is a block coordinate method applied to a nonconvex objective. Such methods generally guarantee descent of the objective under exact block minimization, but not global optimality without additional assumptions. General convergence results for block coordinate descent are due to Grippo and Sciandrone \citep{GrippoSciandrone2000}. Tseng gives convergence theory for broader block coordinate descent settings, including nondifferentiable cases \citep{Tseng2001}.

\subsection{Computational Cost}

For dense \(A\in\mathbb{R}^{m\times n}\), the dominant operations are multiplications of \(A\) or \(A^\top\) by thin matrices. One full alternating iteration costs
\[
O(mnr+(m+n)r^2+r^3).
\]
For a square \(n\times n\) matrix, this becomes
\[
O(n^2r)
\]
when \(r\ll n\). For sparse \(A\), the cost becomes
\[
O(\operatorname{nnz}(A)r+(m+n)r^2+r^3),
\]
where \(\operatorname{nnz}(A)\) is the number of nonzero entries of \(A\).

This is a per-iteration cost. It is not a total runtime guarantee unless the number of iterations is controlled or reported. Iterative methods for sparse linear systems have their own complexity and convergence behavior \citep{Saad2003}. Preconditioning can be essential for large ill-conditioned linear systems \citep{Benzi2002}. Randomized low-rank methods may be faster when the main goal is approximate spectral compression \citep{MartinssonTropp2020}. The proposed method is therefore most relevant when \(r\ll \min\{m,n\}\), when a moderate number of alternating iterations is sufficient, and when explicit regularization or core conditioning is part of the modeling goal.

\subsection{Relation to Linear Systems}

When \(A\in\mathbb{R}^{n\times n}\), the factorization \(A\approx PDQ\) may be useful as a compressed surrogate, denoised approximation, or ingredient in a preconditioner. However, if \(r<n\), the product \(PDQ\) has rank at most \(r\) and cannot directly replace a nonsingular full-rank matrix in the system
\[
Ax=b.
\]
A direct solver based only on \(PDQ\) would therefore require additional structure, such as a diagonal correction, low-rank-plus-sparse model, low-rank-plus-identity model, or iterative refinement. This paper focuses on the approximation and conditioning problem. Any linear-system solver derived from the factorization requires separate analysis.

\subsection{Applicability}

The framework is applicable when a low-rank or approximately low-rank representation is meaningful and when factor regularization or core conditioning is part of the modeling goal. Examples include dense matrices, sparse matrices, covariance matrices, kernel matrices, graph-derived matrices, image matrices, recommender-system matrices, and discretized operator matrices. Recommender matrices motivate weighted or partially observed objectives \citep{HarperKonstan2015}. Image matrices provide standard compression and reconstruction benchmarks \citep{LeCun1998}. Sparse scientific matrices motivate storage-aware and matrix-vector-efficient approximations \citep{DavisHu2011}. The framework also connects naturally to tensor methods through matrix unfoldings and multilinear generalizations \citep{KoldaBader2009}. Tensor-train decompositions provide one important high-dimensional extension of matrix factorization ideas \citep{Oseledets2011}. These connections are not pursued in detail here.

\section{Main Results}
\label{sec:main_results}

We now give the main theoretical results for the core-conditioned tri-factorization. The results are stated for a rectangular matrix \(A\in\mathbb{R}^{m\times n}\), since the square case follows by taking \(m=n\). We distinguish carefully between product-level approximation, factor-level identifiability, conditioning of the core, and convergence of the alternating scheme.

Throughout this section, let
\[
X(P,D,Q)=PDQ,
\qquad
P\in\mathbb{R}^{m\times r},\quad
D\in\mathbb{R}^{r\times r},\quad
Q\in\mathbb{R}^{r\times n},
\]
where \(r\leq \min\{m,n\}\). For full observations, consider
\begin{equation}
\Phi_A(P,D,Q)
=
\frac12\|A-PDQ\|_F^2
+
\frac{\alpha_P}{2}\|P\|_F^2
+
\frac{\alpha_D}{2}\|D\|_F^2
+
\frac{\alpha_Q}{2}\|Q\|_F^2,
\label{eq:main_quadratic_objective}
\end{equation}
where \(\alpha_P,\alpha_D,\alpha_Q>0\). For weighted or partially observed data, replace the residual by
\begin{equation}
\frac12\|W\odot(A-PDQ)\|_F^2,
\label{eq:main_weighted_residual}
\end{equation}
where \(W\in\mathbb{R}^{m\times n}\) is nonnegative. The full-data case is recovered by taking \(W\) to be the all-ones matrix.

\subsection{Rank-\(r\) Benchmark and Regularized Oracle Inequality}

The next result fixes the relationship with truncated SVD. It prevents the method from being interpreted as improving on SVD for the classical unregularized Frobenius problem, where SVD is already optimal by the Eckart--Young--Mirsky theorem \citep{EckartYoung1936,Mirsky1960}.

\begin{theorem}[Rank-\(r\) benchmark and oracle inequality]
\label{thm:rank_benchmark_oracle}
Let \(A\in\mathbb{R}^{m\times n}\), and let \(r\leq \min\{m,n\}\).

\begin{enumerate}[label=(\alph*)]
\item In the unregularized full-data case,
\begin{equation}
\inf_{P,D,Q}\|A-PDQ\|_F^2
=
\sum_{j>r}\sigma_j(A)^2.
\label{eq:svd_benchmark}
\end{equation}
Hence the product \(PDQ\) cannot improve upon truncated SVD for the classical rank-\(r\) Frobenius approximation problem.

\item Let
\begin{equation}
\Psi_A(P,D,Q)
=
\frac12\|W\odot(A-PDQ)\|_F^2
+
\mathcal{R}(P,D,Q),
\label{eq:general_weighted_objective}
\end{equation}
where \(\mathcal{R}\geq0\). If \((P_\star,D_\star,Q_\star)\) is a global minimizer, then for every feasible triple \((P,D,Q)\),
\begin{equation}
\frac12\|W\odot(A-P_\star D_\star Q_\star)\|_F^2
+
\mathcal{R}(P_\star,D_\star,Q_\star)
\leq
\frac12\|W\odot(A-PDQ)\|_F^2
+
\mathcal{R}(P,D,Q).
\label{eq:oracle_inequality}
\end{equation}
\end{enumerate}
\end{theorem}

\begin{proof}
For (a), every product \(PDQ\) has rank at most \(r\). Conversely, every matrix of rank at most \(r\) can be written as \(PDQ\) by using its compact singular value decomposition. Therefore the infimum over triples \((P,D,Q)\) is exactly the infimum over all matrices of rank at most \(r\). Equation \eqref{eq:svd_benchmark} follows from the Eckart--Young--Mirsky theorem.

For (b), the inequality is the defining variational inequality for a global minimizer of \(\Psi_A\).
\end{proof}

\begin{remark}
The useful regime of the proposed method is therefore not the unregularized SVD problem. Its role is in weighted, noisy, structured, regularized, or conditioning-sensitive approximation problems, where the optimization problem differs from the classical SVD objective.
\end{remark}

\subsection{Core Conditioning and Stability of the Block Systems}

A bound on \(\kappa(D)\) does not follow from the factorization alone. To obtain a conditioning guarantee, the admissible core set or the regularization must impose such control explicitly. We use the admissible-set formulation because it makes the assumption transparent.

For \(0<\delta\leq M<\infty\), define
\begin{equation}
\mathcal{D}_{\delta,M}
=
\left\{
D\in\mathbb{R}^{r\times r}:
\sigma_{\min}(D)\geq \delta,\ 
\|D\|_2\leq M
\right\}.
\label{eq:core_conditioned_set}
\end{equation}
The constraint \(\mathcal{D}_{\delta,M}\) is a modeling choice, not a consequence of the factorization.

\begin{theorem}[Core conditioning and block-system bounds]
\label{thm:core_conditioning}
Assume that \(D\in\mathcal{D}_{\delta,M}\), with \(0<\delta\leq M<\infty\), and consider the quadratic objective \eqref{eq:main_quadratic_objective}. Then every feasible triple satisfies
\begin{equation}
\kappa(D)\leq \frac{M}{\delta}.
\label{eq:kappa_D_bound}
\end{equation}
Let \((P_\star,D_\star,Q_\star)\) be a minimizer of the constrained problem and define
\begin{equation}
C_A
=
\frac12\|A\|_F^2+\frac{\alpha_D}{2}r\delta^2.
\label{eq:CA_def}
\end{equation}
Then
\begin{equation}
\|P_\star\|_F^2\leq \frac{2C_A}{\alpha_P},
\qquad
\|Q_\star\|_F^2\leq \frac{2C_A}{\alpha_Q}.
\label{eq:PQ_bounds}
\end{equation}
For the full-data block updates, the normal-equation matrices satisfy
\begin{equation}
\kappa\!\left(D_\star Q_\star Q_\star^\top D_\star^\top+\alpha_P I_r\right)
\leq
1+\frac{2M^2C_A}{\alpha_P\alpha_Q},
\label{eq:P_block_condition}
\end{equation}
and
\begin{equation}
\kappa\!\left(D_\star^\top P_\star^\top P_\star D_\star+\alpha_Q I_r\right)
\leq
1+\frac{2M^2C_A}{\alpha_Q\alpha_P}.
\label{eq:Q_block_condition}
\end{equation}
Moreover, the linear operator
\begin{equation}
\mathcal{L}_{P_\star,Q_\star}(Y)
=
(P_\star^\top P_\star)Y(Q_\star Q_\star^\top)+\alpha_DY
\end{equation}
appearing in the \(D\)-update satisfies
\begin{equation}
\kappa(\mathcal{L}_{P_\star,Q_\star})
\leq
1+
\frac{4C_A^2}{\alpha_D\alpha_P\alpha_Q}.
\label{eq:D_operator_condition}
\end{equation}
\end{theorem}

\begin{proof}
The bound \eqref{eq:kappa_D_bound} follows directly from the definition of \(\mathcal{D}_{\delta,M}\).

The triple \((0,\delta I_r,0)\) is feasible and has objective value \(C_A\). Hence the minimizing property gives
\begin{equation}
\Phi_A(P_\star,D_\star,Q_\star)\leq C_A.
\end{equation}
Since all three quadratic regularization terms are nonnegative,
\begin{equation}
\frac{\alpha_P}{2}\|P_\star\|_F^2\leq C_A,
\qquad
\frac{\alpha_Q}{2}\|Q_\star\|_F^2\leq C_A,
\end{equation}
which proves \eqref{eq:PQ_bounds}.

Let
\begin{equation}
H_P=D_\star Q_\star Q_\star^\top D_\star^\top+\alpha_P I_r.
\end{equation}
Then
\begin{equation}
\lambda_{\min}(H_P)\geq \alpha_P,
\qquad
\lambda_{\max}(H_P)
\leq
\alpha_P+\|D_\star Q_\star\|_2^2
\leq
\alpha_P+M^2\|Q_\star\|_F^2.
\end{equation}
Using \eqref{eq:PQ_bounds} gives \eqref{eq:P_block_condition}. The proof of \eqref{eq:Q_block_condition} is identical.

For the \(D\)-update, vectorization gives
\begin{equation}
\operatorname{vec}(\mathcal{L}_{P_\star,Q_\star}(Y))
=
\left[
(Q_\star Q_\star^\top)\otimes(P_\star^\top P_\star)
+
\alpha_D I_{r^2}
\right]\operatorname{vec}(Y).
\end{equation}
The smallest eigenvalue is at least \(\alpha_D\), while the largest is at most
\begin{equation}
\alpha_D+\|P_\star\|_2^2\|Q_\star\|_2^2
\leq
\alpha_D+\|P_\star\|_F^2\|Q_\star\|_F^2.
\end{equation}
Using \eqref{eq:PQ_bounds} gives \eqref{eq:D_operator_condition}.
\end{proof}

\begin{remark}
The theorem gives a verifiable numerical-stability guarantee for the block systems. It does not claim that the whole decomposition is unique, nor that \(\kappa(D)\) follows from \(\kappa(A)\). The conditioning is controlled because it is built into the admissible core model.
\end{remark}

\subsection{Convergence of Exact Alternating Minimization}

We state the convergence result for the full-data quadratic objective without the core-projection step. The projected or backtracked core-bounded variant is treated algorithmically and empirically in this paper. A separate convergence proof for the constrained projected variant would require additional assumptions on the projection and acceptance rule.

Given \((P^t,D^t,Q^t)\), define one exact alternating-minimization sweep by
\begin{equation}
P^{t+1}
=
\arg\min_P \Phi_A(P,D^t,Q^t),
\label{eq:P_exact_update}
\end{equation}
\begin{equation}
D^{t+1}
=
\arg\min_D \Phi_A(P^{t+1},D,Q^t),
\label{eq:D_exact_update}
\end{equation}
and
\begin{equation}
Q^{t+1}
=
\arg\min_Q \Phi_A(P^{t+1},D^{t+1},Q).
\label{eq:Q_exact_update}
\end{equation}
For \(\alpha_P,\alpha_D,\alpha_Q>0\), each block subproblem is strongly convex and has a unique minimizer.

\begin{theorem}[Descent and critical-point convergence for the full-data quadratic case]
\label{thm:alternating_convergence}
Let \(A\in\mathbb{R}^{m\times n}\), and let \(\alpha_P,\alpha_D,\alpha_Q>0\). Let \(\{(P^t,D^t,Q^t)\}_{t\geq0}\) be generated by the exact alternating scheme \eqref{eq:P_exact_update}--\eqref{eq:Q_exact_update} for the objective \eqref{eq:main_quadratic_objective}. Then the following hold.

\begin{enumerate}[label=(\alph*)]
\item The objective values are nonincreasing:
\[
\Phi_A(P^{t+1},D^{t+1},Q^{t+1})
\leq
\Phi_A(P^t,D^t,Q^t).
\]

\item The sequence \(\{(P^t,D^t,Q^t)\}_{t\geq0}\) is bounded.

\item The squared successive differences are summable:
\[
\sum_{t=0}^{\infty}
\left(
\|P^{t+1}-P^t\|_F^2
+
\|D^{t+1}-D^t\|_F^2
+
\|Q^{t+1}-Q^t\|_F^2
\right)
<\infty.
\]

\item Every cluster point of the sequence is a critical point of \(\Phi_A\).

\item Since \(\Phi_A\) is a polynomial function of the entries of \(P,D,Q\), it satisfies the Kurdyka--\L{}ojasiewicz property. Consequently, the whole sequence converges to a critical point of \(\Phi_A\).
\end{enumerate}
\end{theorem}

\begin{proof}
Each block subproblem is strongly convex because the updated block appears with a positive quadratic regularization term. Exact minimization gives
\[
\Phi_A(P^t,D^t,Q^t)
\geq
\Phi_A(P^{t+1},D^t,Q^t)
\geq
\Phi_A(P^{t+1},D^{t+1},Q^t)
\geq
\Phi_A(P^{t+1},D^{t+1},Q^{t+1}),
\]
which proves monotonicity.

Since \(\Phi_A(P^t,D^t,Q^t)\leq \Phi_A(P^0,D^0,Q^0)\) and
\[
\Phi_A(P,D,Q)
\geq
\frac{\alpha_P}{2}\|P\|_F^2
+
\frac{\alpha_D}{2}\|D\|_F^2
+
\frac{\alpha_Q}{2}\|Q\|_F^2,
\]
the iterates remain in a bounded level set.

Strong convexity of each block subproblem gives sufficient decrease:
\[
\Phi_A(P^t,D^t,Q^t)-\Phi_A(P^{t+1},D^{t+1},Q^{t+1})
\geq
c
\left(
\|P^{t+1}-P^t\|_F^2
+
\|D^{t+1}-D^t\|_F^2
+
\|Q^{t+1}-Q^t\|_F^2
\right),
\]
where \(c=\frac12\min\{\alpha_P,\alpha_D,\alpha_Q\}\). Summing this inequality over \(t\) proves the summability of squared successive differences.

The gradient of \(\Phi_A\) is Lipschitz on the bounded level set. The first-order optimality conditions of the exact block updates imply a relative-error estimate of the form
\[
\|\nabla\Phi_A(P^{t+1},D^{t+1},Q^{t+1})\|
\leq
C
\left(
\|P^{t+1}-P^t\|_F
+
\|D^{t+1}-D^t\|_F
+
\|Q^{t+1}-Q^t\|_F
\right)
\]
for some constant \(C>0\) depending on the bounded level set. Since the right-hand side tends to zero, every cluster point is critical.

The function \(\Phi_A\) is polynomial in the entries of \(P,D,Q\), and therefore satisfies the Kurdyka--\L{}ojasiewicz property. The sufficient-decrease estimate, the relative-error estimate, boundedness, and the KL property imply finite length and convergence of the whole sequence to a critical point by the standard KL convergence theorem for descent methods \citep{AttouchBolteSvaiter2013,BolteSabachTeboulle2014}.
\end{proof}

\begin{remark}
The theorem proves convergence to a critical point for the exact full-data quadratic alternating scheme. It does not prove convergence to a global minimizer. It also does not cover the weighted missing-entry implementation or the projected core-bounded implementation without additional assumptions.
\end{remark}

\subsection{Perturbation Stability}

The next result gives a product-level perturbation bound. This is the right stability statement without additional gauge fixing. A factor-level bound such as \(\|D-D_0\|_F\leq C\varepsilon\) is generally false unless the factorization has been made locally identifiable.

\begin{theorem}[Product-level stability under perturbations]
\label{thm:product_stability}
Let \(A=A_0+E\), and suppose
\[
\|W\odot E\|_F\leq \varepsilon.
\]
Let
\[
\Psi_A(P,D,Q)
=
\frac12\|W\odot(A-PDQ)\|_F^2
+
\mathcal{R}(P,D,Q),
\]
where \(\mathcal{R}\geq0\). Let \((P_\star,D_\star,Q_\star)\) be a global minimizer of \(\Psi_A\). For any feasible reference triple \((P_0,D_0,Q_0)\), set
\[
X_\star=P_\star D_\star Q_\star,
\qquad
X_0=P_0D_0Q_0,
\qquad
\eta_0=\|W\odot(A_0-X_0)\|_F.
\]
Then
\begin{equation}
\|W\odot(A_0-X_\star)\|_F
\leq
\varepsilon
+
\sqrt{(\eta_0+\varepsilon)^2+2\mathcal{R}(P_0,D_0,Q_0)}.
\label{eq:product_stability_bound}
\end{equation}
In particular, if \(A_0=X_0\) on the observed entries, then
\begin{equation}
\|W\odot(A_0-X_\star)\|_F
\leq
\varepsilon
+
\sqrt{\varepsilon^2+2\mathcal{R}(P_0,D_0,Q_0)}.
\label{eq:exact_reference_stability}
\end{equation}
\end{theorem}

\begin{proof}
By optimality of \((P_\star,D_\star,Q_\star)\),
\[
\frac12\|W\odot(A-X_\star)\|_F^2
+
\mathcal{R}(P_\star,D_\star,Q_\star)
\leq
\frac12\|W\odot(A-X_0)\|_F^2
+
\mathcal{R}(P_0,D_0,Q_0).
\]
Since \(\mathcal{R}\geq0\),
\[
\|W\odot(A-X_\star)\|_F
\leq
\sqrt{\|W\odot(A-X_0)\|_F^2+2\mathcal{R}(P_0,D_0,Q_0)}.
\]
Also,
\[
\|W\odot(A-X_0)\|_F
=
\|W\odot(A_0-X_0)+W\odot E\|_F
\leq
\eta_0+\varepsilon.
\]
Finally,
\[
\|W\odot(A_0-X_\star)\|_F
\leq
\|W\odot(A-X_\star)\|_F+\varepsilon.
\]
Combining the last three inequalities gives \eqref{eq:product_stability_bound}. The special case \eqref{eq:exact_reference_stability} follows by taking \(\eta_0=0\).
\end{proof}

\begin{corollary}[Local factor stability under gauge fixing]
\label{cor:local_factor_stability}
Let \(\mathcal{G}\) be a smooth gauge-fixed parameter space that removes the invariance
\[
PDQ=(PS)(S^{-1}DT^{-1})(TQ).
\]
Let \(z=(P,D,Q)\in\mathcal{G}\), and assume that \(z_0\) is an isolated local minimizer of \(\Psi_{A_0}\). Suppose that, in a neighborhood \(U\subset\mathcal{G}\) of \(z_0\),
\begin{equation}
\nabla^2\Psi_{A_0}(z)\succeq \mu I
\qquad\text{for some }\mu>0,
\label{eq:local_strong_growth}
\end{equation}
and that the data perturbation satisfies
\begin{equation}
\|\nabla\Psi_A(z)-\nabla\Psi_{A_0}(z)\|
\leq
L\|W\odot(A-A_0)\|_F
\qquad
\text{for all } z\in U.
\label{eq:gradient_perturbation_bound}
\end{equation}
If \(z_A\in U\) is a local minimizer of \(\Psi_A\), then
\begin{equation}
\|z_A-z_0\|
\leq
\frac{L}{\mu}\|W\odot(A-A_0)\|_F.
\label{eq:local_factor_lipschitz}
\end{equation}
Consequently, under these explicit local identifiability assumptions,
\begin{equation}
\|D_A-D_0\|_F
\leq
\frac{L}{\mu}\|W\odot(A-A_0)\|_F.
\label{eq:D_local_stability}
\end{equation}
\end{corollary}

\begin{proof}
The first-order conditions give
\[
\nabla\Psi_A(z_A)=0,
\qquad
\nabla\Psi_{A_0}(z_0)=0.
\]
By \eqref{eq:local_strong_growth},
\[
\mu\|z_A-z_0\|^2
\leq
\left\langle
\nabla\Psi_{A_0}(z_A)-\nabla\Psi_{A_0}(z_0),
z_A-z_0
\right\rangle.
\]
Using the first-order conditions,
\[
\nabla\Psi_{A_0}(z_A)-\nabla\Psi_{A_0}(z_0)
=
\nabla\Psi_{A_0}(z_A)-\nabla\Psi_A(z_A).
\]
Therefore, by \eqref{eq:gradient_perturbation_bound},
\[
\mu\|z_A-z_0\|^2
\leq
L\|W\odot(A-A_0)\|_F\|z_A-z_0\|.
\]
If \(z_A\neq z_0\), divide by \(\|z_A-z_0\|\). The case \(z_A=z_0\) is immediate. The bound for \(D\) follows because \(D\) is one component of \(z\).
\end{proof}

\begin{remark}
This corollary replaces the unsupported global claim \(\|D-D_0\|_F\leq C\varepsilon\). A factor-level perturbation bound is available only after the latent-space symmetry has been removed and a local second-order stability condition has been imposed. Without such assumptions, only product-level stability is generally meaningful.
\end{remark}

\subsection{Computational Scaling}

For dense \(A\in\mathbb{R}^{m\times n}\), one exact full-data sweep of the quadratic alternating scheme requires
\begin{equation}
O(mnr+(m+n)r^2+r^3)
\label{eq:dense_cost}
\end{equation}
floating-point operations. For sparse \(A\),
\begin{equation}
O(\operatorname{nnz}(A)r+(m+n)r^2+r^3)
\label{eq:sparse_cost}
\end{equation}
operations are required, where \(\operatorname{nnz}(A)\) is the number of nonzero entries. These are per-iteration costs. They should be compared with classical and randomized methods according to the actual task. Direct dense factorizations remain appropriate for exact full-rank linear algebra \citep{GolubVanLoan2013}. Randomized methods can be faster for approximate spectral compression \citep{HalkoMartinssonTropp2011,MartinssonTropp2020}. The advantage of the present formulation is the joint control of approximation, factor regularity, and core conditioning.

\section{Algorithmic Framework and Illustrative Examples}
\label{sec:algorithm_examples}

We now give the computational form of the core-conditioned tri-factorization. The algorithm is stated for the quadratic regularized objective
\begin{equation}
\Phi_A(P,D,Q)
=
\frac12\|A-PDQ\|_F^2
+
\frac{\alpha_P}{2}\|P\|_F^2
+
\frac{\alpha_D}{2}\|D\|_F^2
+
\frac{\alpha_Q}{2}\|Q\|_F^2,
\label{eq:algorithm_objective}
\end{equation}
where \(\alpha_P,\alpha_D,\alpha_Q>0\). For missing or weighted data, the residual term is replaced by
\[
\frac12\|W\odot(A-PDQ)\|_F^2.
\]
The full-data version admits closed-form block updates. The weighted version uses the same block-coordinate principle, but each block update becomes a weighted ridge least-squares problem and is treated here as an implementation extension. The convergence theorem in Section~\ref{sec:main_results} applies only to the exact full-data quadratic scheme.

\subsection{Core-Conditioned Alternating Tri-Factorization}

Let \(r\leq \min\{m,n\}\) be the target rank. Given \(D\) and \(Q\), the \(P\)-update in the full-data case is
\begin{equation}
P
=
A Q^\top D^\top
\left(DQQ^\top D^\top+\alpha_P I_r\right)^{-1}.
\label{eq:algorithm_P_update}
\end{equation}
Given \(P\) and \(Q\), the \(D\)-update solves the Sylvester-type equation
\begin{equation}
(P^\top P)D(QQ^\top)+\alpha_DD
=
P^\top A Q^\top.
\label{eq:algorithm_D_update}
\end{equation}
Given \(P\) and \(D\), the \(Q\)-update solves
\begin{equation}
\left(D^\top P^\top P D+\alpha_Q I_r\right)Q
=
D^\top P^\top A.
\label{eq:algorithm_Q_update}
\end{equation}
The positive ridge parameters make the block subproblems well posed. This is a standard block-coordinate structure, so the algorithm should be assessed using objective descent, stationarity, reconstruction error, and conditioning rather than only final Frobenius error \citep{GrippoSciandrone2000,Tseng2001}.

\begin{algorithm}[!htbp]
\caption{Core-Conditioned Alternating Tri-Factorization}
\label{alg:ccd_tf}
\begin{algorithmic}[1]
\Require Matrix \(A\in\mathbb{R}^{m\times n}\), target rank \(r\), parameters \(\alpha_P,\alpha_D,\alpha_Q>0\), tolerance \(\tau>0\), maximum iterations \(T_{\max}\), optional core bounds \(0<\delta\leq M<\infty\).
\Ensure Factors \(P\in\mathbb{R}^{m\times r}\), \(D\in\mathbb{R}^{r\times r}\), \(Q\in\mathbb{R}^{r\times n}\), and diagnostics.

\State Initialize \(P^0,Q^0\) by randomized range approximation or truncated SVD warm start.
\State Initialize \(D^0=I_r\), or \(D^0=\operatorname{diag}(\sigma_1,\ldots,\sigma_r)\) if an SVD warm start is used.
\State Set \(t=0\). Compute \(\Phi_A(P^0,D^0,Q^0)\).

\Repeat
\State Update
\[
P^{t+1}
=
A(Q^t)^\top(D^t)^\top
\left(D^tQ^t(Q^t)^\top(D^t)^\top+\alpha_P I_r\right)^{-1}.
\]

\State Compute \(\widehat D^{t+1}\) as the solution of
\[
((P^{t+1})^\top P^{t+1})D(Q^t(Q^t)^\top)+\alpha_DD
=
(P^{t+1})^\top A(Q^t)^\top.
\]

\If{core bounds \(\mathcal{D}_{\delta,M}\) are imposed}
\State Compute the singular value decomposition \(\widehat D^{t+1}=U\Sigma V^\top\).
\State Set
\[
D^{t+1}
=
U\,\operatorname{diag}\!\left(
\min\{M,\max\{\delta,\sigma_i(\widehat D^{t+1})\}\}
\right)V^\top.
\]
\State Accept the projected update only if it does not increase the objective; otherwise use a backtracking step between \(D^t\) and the projected candidate.
\Else
\State Set \(D^{t+1}=\widehat D^{t+1}\).
\EndIf

\State Update \(Q^{t+1}\) by solving
\[
\left((D^{t+1})^\top(P^{t+1})^\top P^{t+1}D^{t+1}+\alpha_Q I_r\right)Q
=
(D^{t+1})^\top(P^{t+1})^\top A.
\]

\State Compute the relative reconstruction error
\[
\operatorname{relerr}_{t+1}
=
\frac{\|A-P^{t+1}D^{t+1}Q^{t+1}\|_F}{\|A\|_F}.
\]

\State Compute the relative objective decrease
\[
\Delta_{t+1}
=
\frac{\Phi_A(P^t,D^t,Q^t)-\Phi_A(P^{t+1},D^{t+1},Q^{t+1})}
{\max\{1,\Phi_A(P^t,D^t,Q^t)\}}.
\]

\State Compute the relative step size
\[
s_{t+1}
=
\frac{
\|P^{t+1}-P^t\|_F+\|D^{t+1}-D^t\|_F+\|Q^{t+1}-Q^t\|_F
}
{\max\{1,\|P^t\|_F+\|D^t\|_F+\|Q^t\|_F\}}.
\]

\State Compute the core condition number
\[
\kappa_D^{t+1}=\kappa(D^{t+1}).
\]

\State Set \(t\leftarrow t+1\).
\Until{\(\max\{\Delta_t,s_t\}\leq \tau\) or \(t=T_{\max}\)}

\State \Return \(P^t,D^t,Q^t\), \(\operatorname{relerr}_t\), \(\Phi_A(P^t,D^t,Q^t)\), \(\kappa_D^t\), \(s_t\), \(\Delta_t\).
\end{algorithmic}
\end{algorithm}

\subsection{Diagnostics and Complexity}

The algorithm should be reported with the following diagnostics:
\[
\operatorname{relerr}
=
\frac{\|A-PDQ\|_F}{\|A\|_F},
\qquad
\operatorname{obj}
=
\Phi_A(P,D,Q),
\]
\[
\kappa_D=\kappa(D),
\qquad
s_t=
\frac{\|P^{t}-P^{t-1}\|_F+\|D^{t}-D^{t-1}\|_F+\|Q^{t}-Q^{t-1}\|_F}
{\max\{1,\|P^{t-1}\|_F+\|D^{t-1}\|_F+\|Q^{t-1}\|_F\}},
\]
and
\[
\Delta_t=
\frac{\Phi_A(P^{t-1},D^{t-1},Q^{t-1})-\Phi_A(P^t,D^t,Q^t)}
{\max\{1,\Phi_A(P^{t-1},D^{t-1},Q^{t-1})\}}.
\]
For dense \(A\in\mathbb{R}^{m\times n}\), one full sweep costs
\[
O(mnr+(m+n)r^2+r^3).
\]
For sparse \(A\), the cost becomes
\[
O(\operatorname{nnz}(A)r+(m+n)r^2+r^3).
\]
These are per-sweep costs, not total-runtime guarantees. Total runtime depends on the number of sweeps required to reach the stopping tolerance.

\section{Illustrative Examples}
\label{sec:examples}

The following examples are not intended to replace numerical experiments. Their role is to clarify what the factorization represents, what the conditioning mechanism does, and which claims are justified.

\begin{example}[Exact representation at sufficient rank]
\label{ex:exact_rank}
Let
\[
A=
\begin{bmatrix}
3 & 0 & 3\\
0 & 1 & 1\\
3 & 1 & 4
\end{bmatrix}.
\]
This matrix has rank \(2\), since
\[
A=
U\Sigma U^\top,\qquad
U=
\begin{bmatrix}
1 & 0\\
0 & 1\\
1 & 1
\end{bmatrix},
\qquad
\Sigma=
\begin{bmatrix}
3 & 0\\
0 & 1
\end{bmatrix}.
\]
Taking
\[
P=U,\qquad D=\Sigma,\qquad Q=U^\top
\]
gives
\[
PDQ=U\Sigma U^\top=A.
\]
Thus, when the target rank is at least the true rank, exact representation is possible. This example also shows the correct interpretation of the model: \(PDQ\) is a rank-constrained representation. Its value comes from regularization and conditioning, not from enlarging the rank-\(r\) algebraic class.
\end{example}

\begin{example}[Core conditioning through gauge choice]
\label{ex:core_conditioning}
Let
\[
A=
\begin{bmatrix}
1 & 0\\
0 & \varepsilon
\end{bmatrix},
\qquad
0<\varepsilon\ll1.
\]
One exact representation is
\[
P=I_2,\qquad
D=
\begin{bmatrix}
1 & 0\\
0 & \varepsilon
\end{bmatrix},
\qquad
Q=I_2.
\]
This gives \(PDQ=A\), but
\[
\kappa(D)=\varepsilon^{-1}.
\]
A second exact representation is
\[
P=
\begin{bmatrix}
1 & 0\\
0 & \sqrt{\varepsilon}
\end{bmatrix},
\qquad
D=I_2,
\qquad
Q=
\begin{bmatrix}
1 & 0\\
0 & \sqrt{\varepsilon}
\end{bmatrix}.
\]
Again \(PDQ=A\), but now
\[
\kappa(D)=1.
\]
The same product can therefore be represented with either an ill-conditioned or a perfectly conditioned core. This illustrates why core conditioning must be imposed deliberately. It is not inherited automatically from the matrix \(A\), and it is not guaranteed by the factorization alone.
\end{example}

\begin{example}[Weighted approximation with missing observations]
\label{ex:weighted_missing}
Let
\[
A=
\begin{bmatrix}
5 & ? & 1\\
4 & 3 & ?\\
? & 2 & 1
\end{bmatrix},
\qquad
W=
\begin{bmatrix}
1 & 0 & 1\\
1 & 1 & 0\\
0 & 1 & 1
\end{bmatrix}.
\]
The symbol \(?\) denotes an unobserved entry. The appropriate objective is not the full Frobenius loss, but
\[
\min_{P,D,Q}
\frac12\|W\odot(A-PDQ)\|_F^2
+
\mathcal{R}(P,D,Q).
\]
The observed-entry residual is
\[
\|W\odot(A-PDQ)\|_F^2
=
\sum_{(i,j):W_{ij}=1}(A_{ij}-(PDQ)_{ij})^2.
\]
Thus the method fits only the measured entries while using the low-rank structure and regularization to estimate the missing ones. This is the correct formulation for matrix completion or recommender-type data. It is different from truncated SVD on a fully observed matrix, and the comparison should be made against methods designed for missing-data low-rank approximation.
\end{example}

\begin{example}[Perturbation-stable product recovery]
\label{ex:perturbation}
Let
\[
A_0=P_0D_0Q_0,
\qquad
\|E\|_F\leq\varepsilon,
\qquad
A=A_0+E.
\]
Assume full observations and let \((P_\star,D_\star,Q_\star)\) minimize
\[
\Psi_A(P,D,Q)
=
\frac12\|A-PDQ\|_F^2+\mathcal{R}(P,D,Q),
\qquad
\mathcal{R}\geq0.
\]
Taking \(X_\star=P_\star D_\star Q_\star\) and \(X_0=P_0D_0Q_0=A_0\), Theorem~\ref{thm:product_stability} gives
\[
\|A_0-X_\star\|_F
\leq
\varepsilon+\sqrt{\varepsilon^2+2\mathcal{R}(P_0,D_0,Q_0)}.
\]
If the reference factors are lightly regularized and \(\mathcal{R}(P_0,D_0,Q_0)=O(\varepsilon^2)\), then
\[
\|A_0-X_\star\|_F=O(\varepsilon).
\]
This is the correct perturbation claim at the product level. A bound on the individual core \(D\) requires additional gauge fixing and local identifiability assumptions.
\end{example}

\begin{example}[Dense and sparse computational scaling]
\label{ex:scaling}
Let \(A\in\mathbb{R}^{3000\times2000}\) and choose \(r=100\). For dense \(A\), one full alternating sweep costs
\[
O(mnr+(m+n)r^2+r^3)
=
O(3000\cdot2000\cdot100+5000\cdot100^2+100^3).
\]
The leading term is
\[
6.0\times10^8
\]
floating-point-scale operations. If instead \(A\) is sparse with \(\operatorname{nnz}(A)=2.0\times10^6\), then one sweep costs
\[
O(\operatorname{nnz}(A)r+(m+n)r^2+r^3)
=
O(2.0\times10^6\cdot100+5000\cdot100^2+100^3),
\]
whose leading term is
\[
2.0\times10^8.
\]
These are per-sweep costs. The comparison with LU, QR, SVD, randomized SVD, or Krylov methods depends on the task. The advantage of the present method is not that it universally beats classical factorizations, but that it provides a controllable low-rank representation with explicit regularization and core-conditioning diagnostics.
\end{example}

\section{Numerical Validation}
\label{sec:numerical_validation}

This section evaluates the proposed core-conditioned tri-factorization in controlled matrix regimes where low-rank approximation, regularization, and core conditioning can be assessed directly. The experiments cover clean low-rank recovery, noisy low-rank recovery, ill-conditioned spectra, missing-entry recovery, ablation of the regularization structure, and runtime scaling. All reported values are from the full multi-seed validation run. For each experiment, values are reported as mean \(\pm\) standard deviation across repeated seeds when repetitions are available.

The goal is not to show that the method improves on truncated SVD for the classical full-data unregularized Frobenius rank-\(r\) problem. In that setting, truncated SVD is the correct optimal benchmark. The goal is narrower: to test whether the proposed regularized tri-factorization gives competitive approximation while exposing additional diagnostics, especially the learned core condition number \(\kappa(D)\), convergence behavior, and the effect of side-factor and core regularization.

For full-data experiments, we report
\[
\operatorname{RelFrob}
=
\frac{\|A-\widehat A\|_F}{\|A\|_F},
\qquad
\operatorname{CleanErr}
=
\frac{\|A_0-\widehat A\|_F}{\|A_0\|_F},
\]
where \(A_0\) denotes the noiseless target when available and \(\widehat A=PDQ\). For missing-entry experiments, train, validation, and test RMSE are computed on their corresponding observed-entry splits. Runtime is wall-clock time in seconds, and sweeps denote the number of alternating updates.

The main text reports aggregate performance figures and one representative set of convergence histories using the first prespecified seed. The remaining per-seed histories are placed in Appendix~\ref{app:numerical_evidence}. This keeps the main text readable while preserving the audit trail of the multi-seed run.

\begin{figure}[!htbp]
\centering
\begin{subfigure}[t]{0.24\textwidth}
\includegraphics[width=\linewidth]{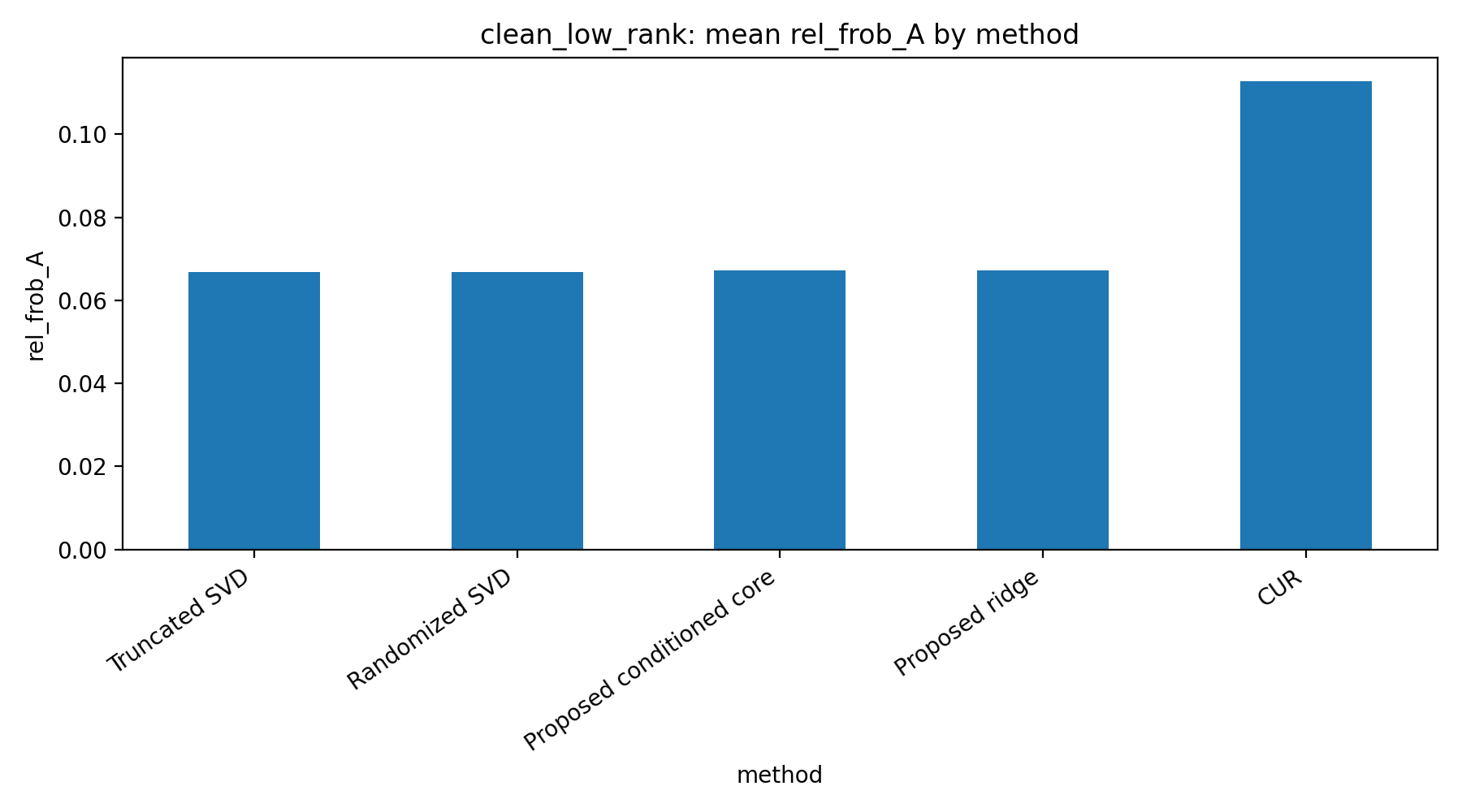}
\caption{Clean low-rank recovery.}
\label{fig:main_clean_low_rank}
\end{subfigure}\hfill
\begin{subfigure}[t]{0.24\textwidth}
\includegraphics[width=\linewidth]{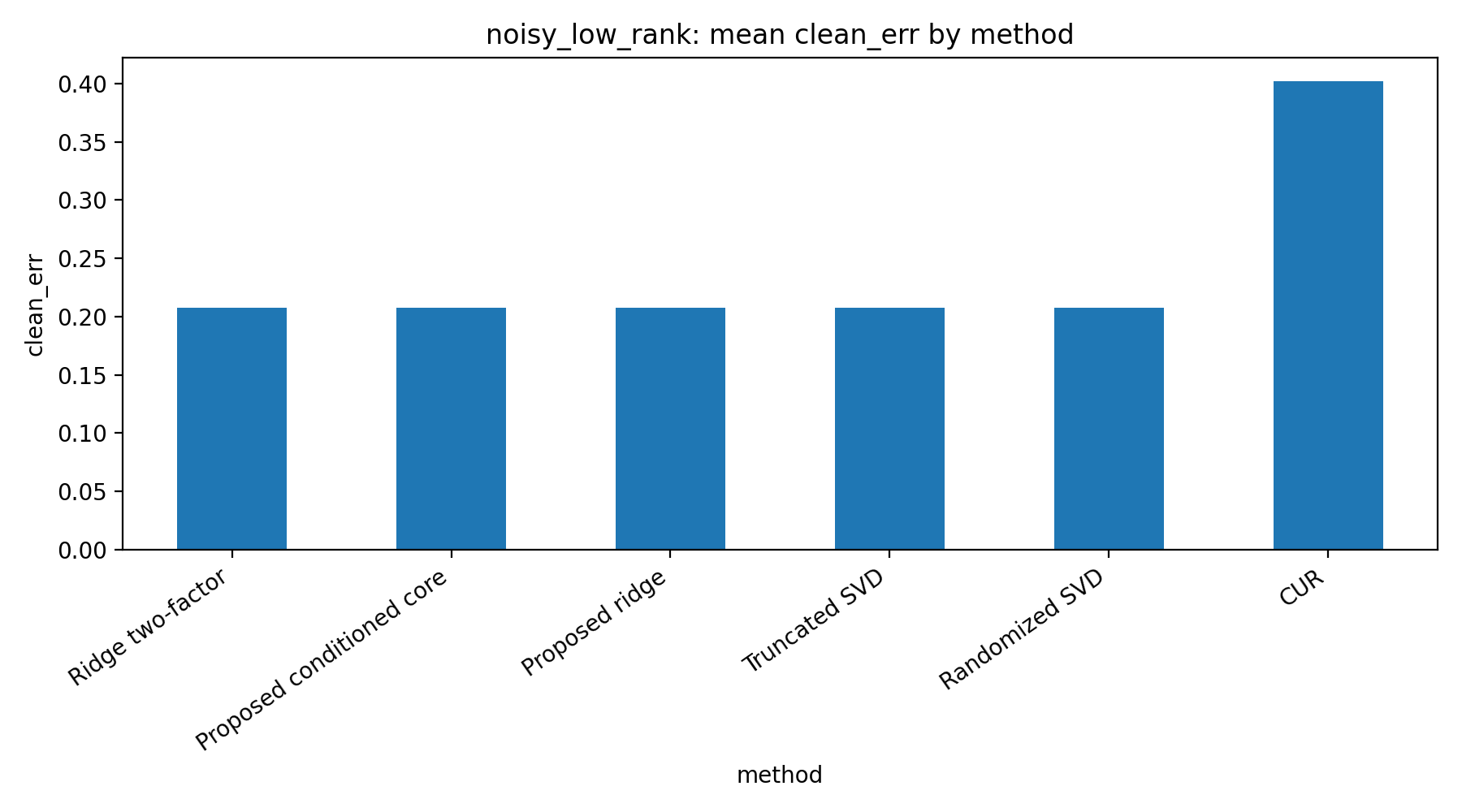}
\caption{Noisy low-rank recovery.}
\label{fig:main_noisy_low_rank}
\end{subfigure}\hfill
\begin{subfigure}[t]{0.24\textwidth}
\includegraphics[width=\linewidth]{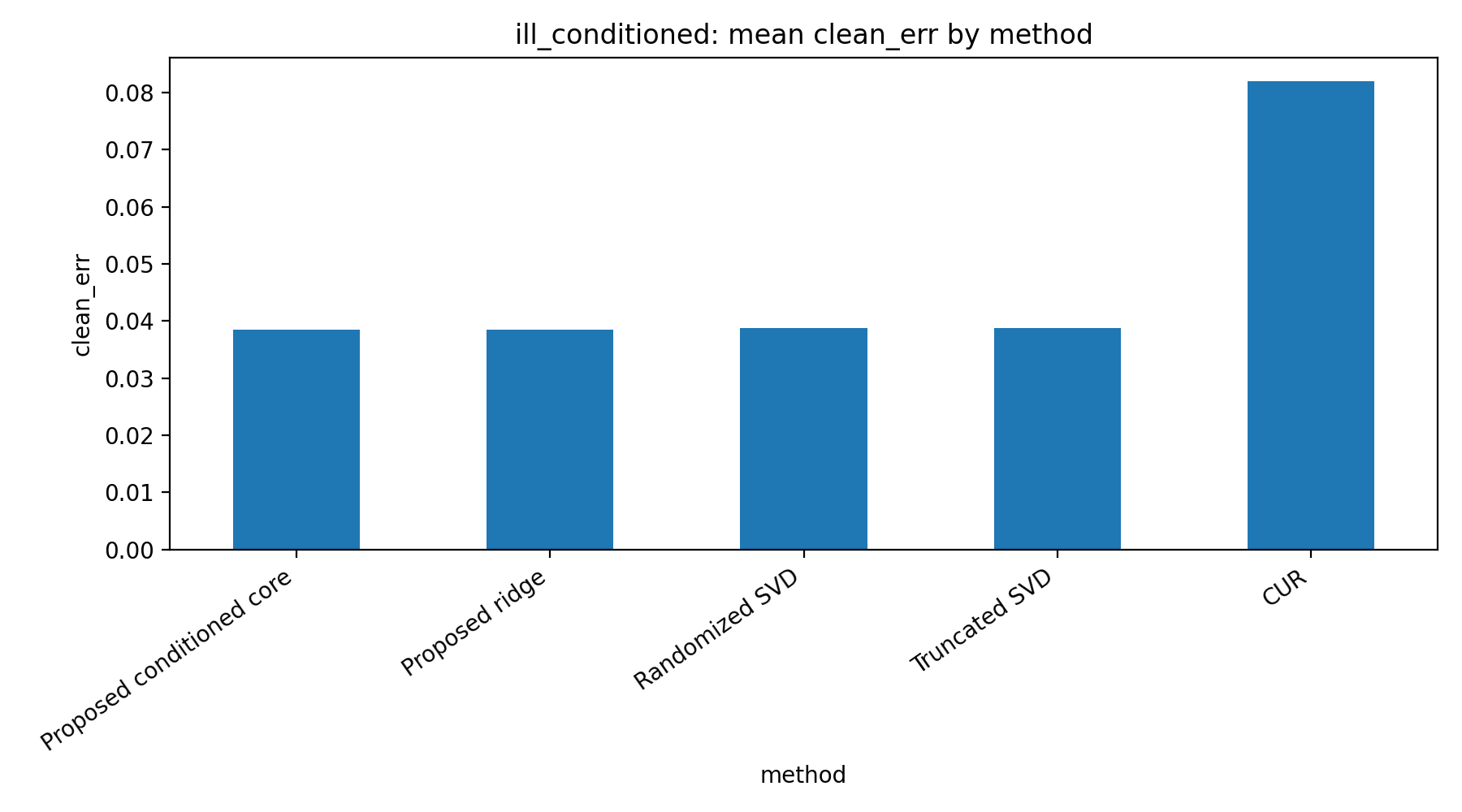}
\caption{Ill-conditioned recovery.}
\label{fig:main_ill_conditioned_error}
\end{subfigure}\hfill
\begin{subfigure}[t]{0.24\textwidth}
\includegraphics[width=\linewidth]{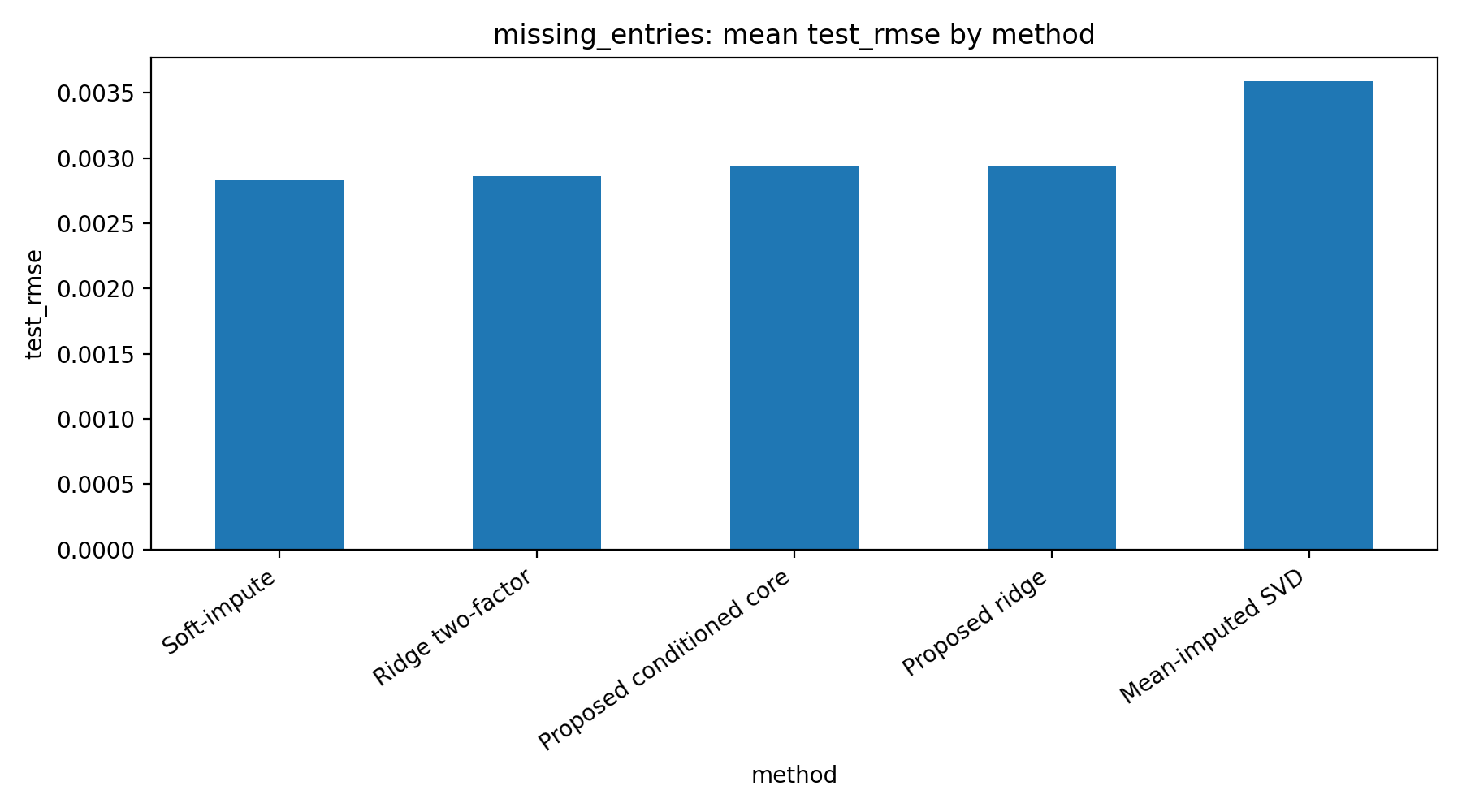}
\caption{Missing-entry test RMSE.}
\label{fig:main_missing_entries}
\end{subfigure}

\medskip

\begin{subfigure}[t]{0.24\textwidth}
\includegraphics[width=\linewidth]{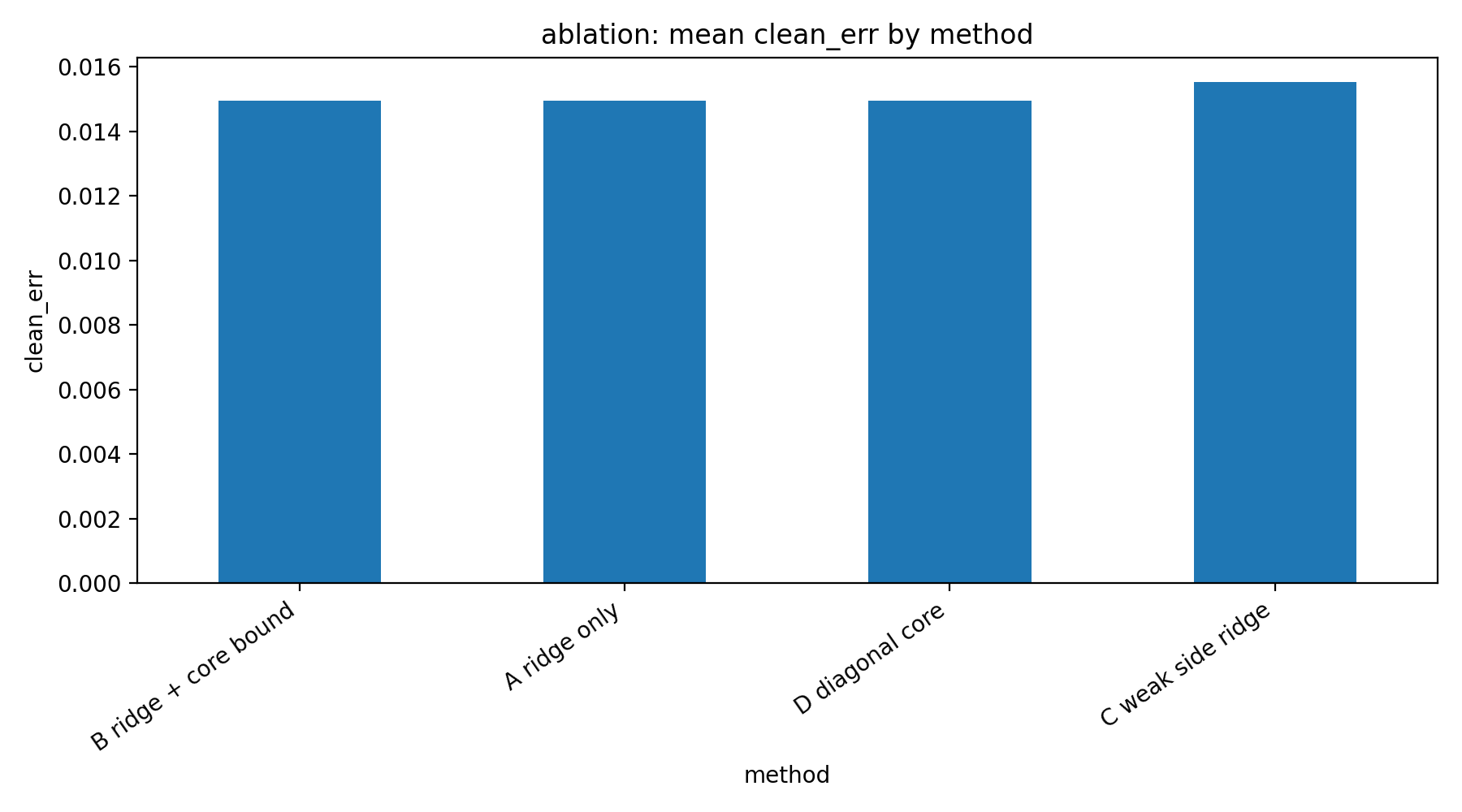}
\caption{Ablation clean error.}
\label{fig:main_ablation_clean_error}
\end{subfigure}\hfill
\begin{subfigure}[t]{0.24\textwidth}
\includegraphics[width=\linewidth]{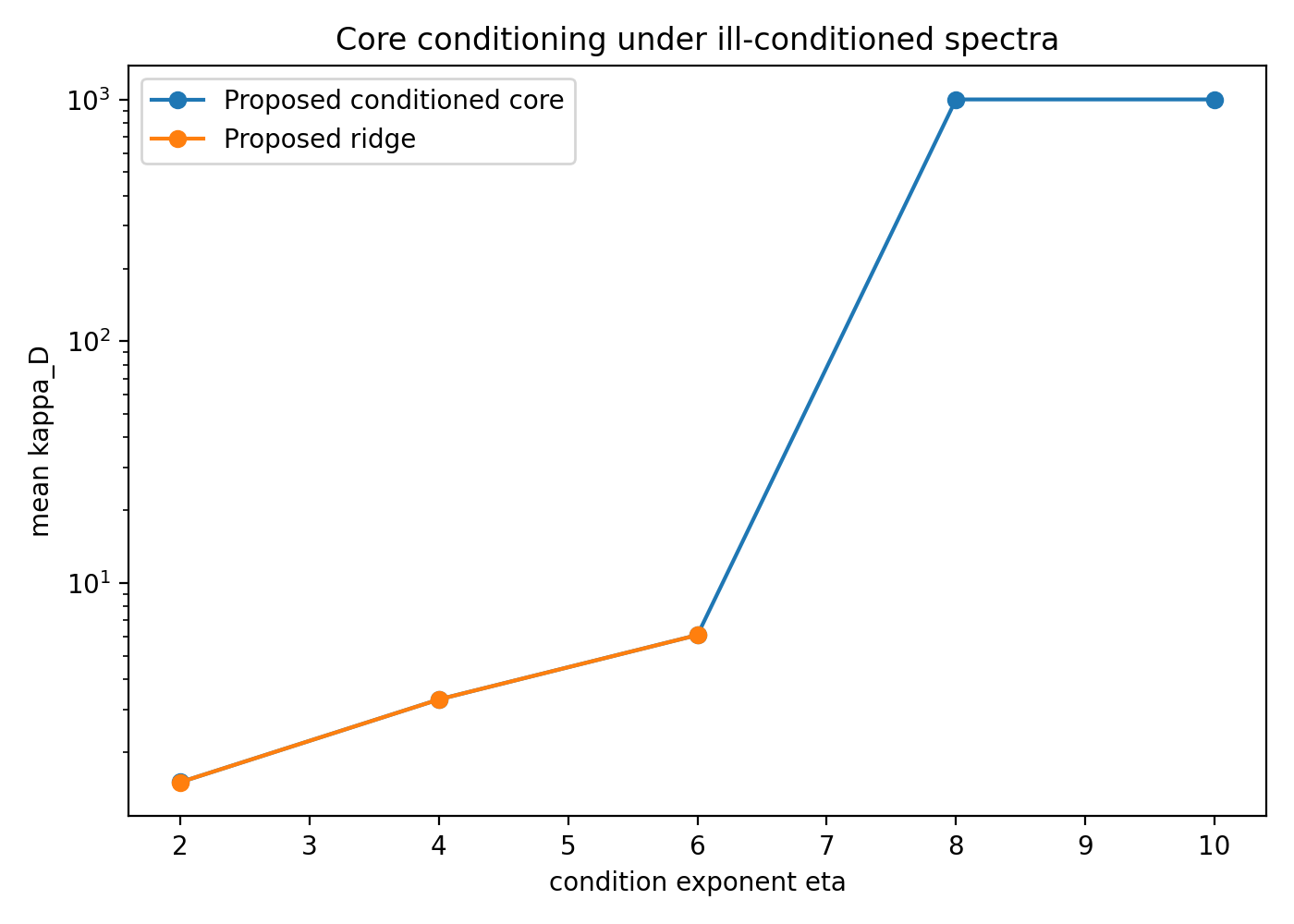}
\caption{Core condition number.}
\label{fig:main_ill_conditioned_kappaD}
\end{subfigure}\hfill
\begin{subfigure}[t]{0.24\textwidth}
\includegraphics[width=\linewidth]{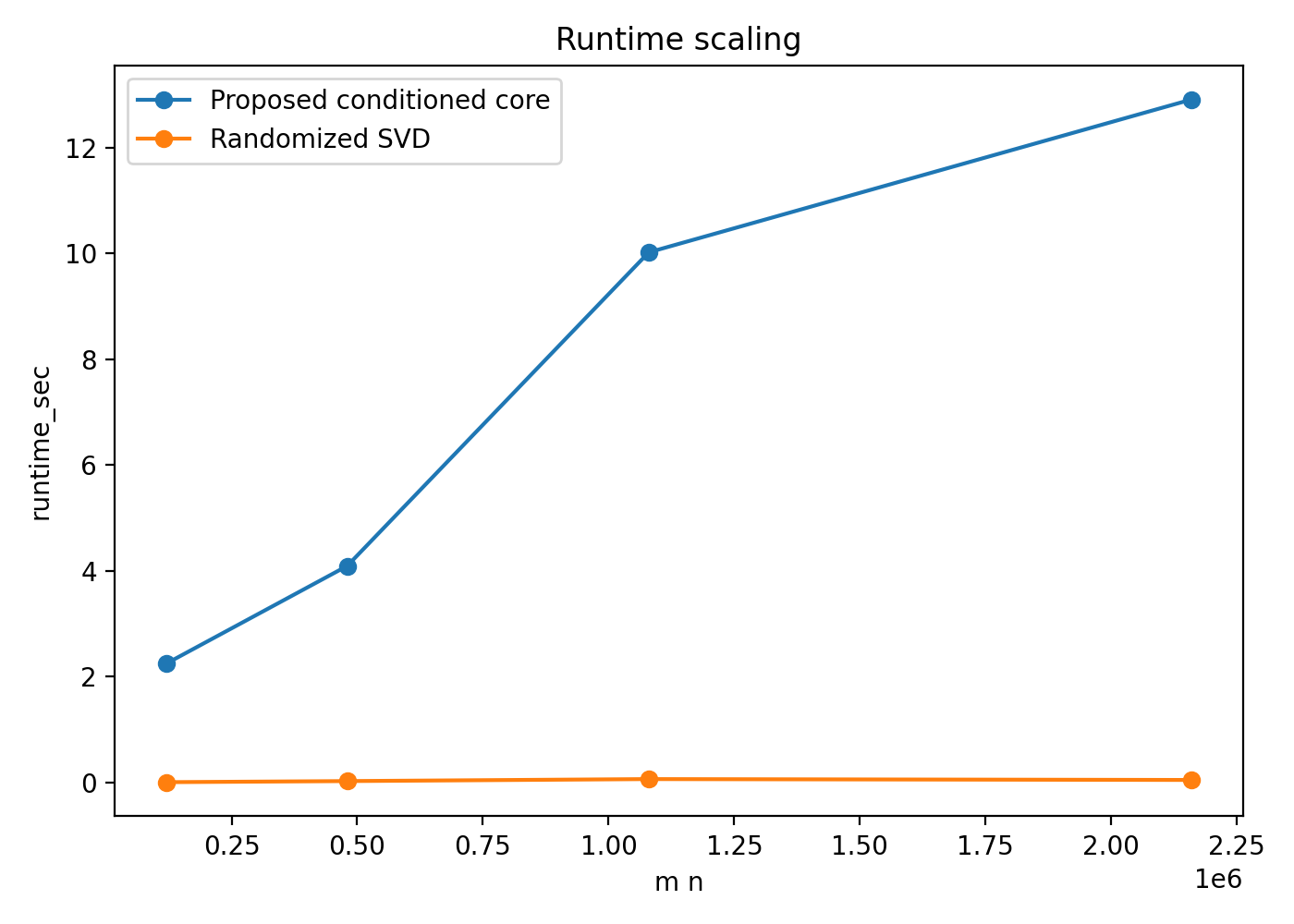}
\caption{Runtime scaling.}
\label{fig:main_runtime_scaling}
\end{subfigure}\hfill
\begin{subfigure}[t]{0.24\textwidth}
\includegraphics[width=\linewidth]{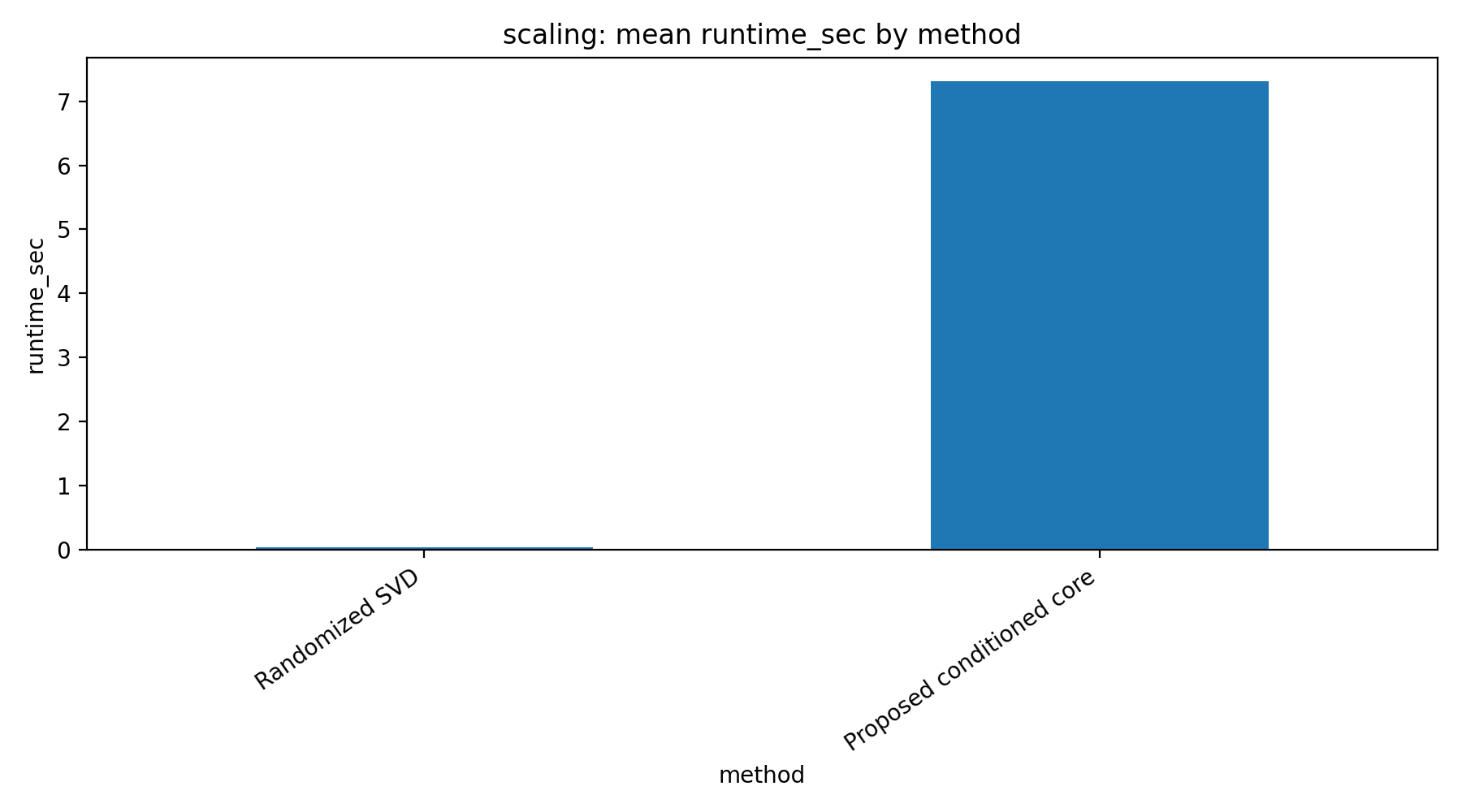}
\caption{Runtime by method.}
\label{fig:main_scaling_runtime}
\end{subfigure}

\caption{Aggregate numerical validation over the full multi-seed run. Panels (a)--(e) summarize reconstruction, noisy recovery, ill-conditioned recovery, missing-entry prediction, and ablation behavior. Panels (f)--(h) summarize core conditioning and runtime behavior. These plots provide the main empirical overview; detailed per-seed histories are reported in Appendix~\ref{app:numerical_evidence}.}
\label{fig:main_aggregate_validation}
\end{figure}

Figure~\ref{fig:main_aggregate_validation} gives the aggregate view of the validation. Panels~\ref{fig:main_clean_low_rank}--\ref{fig:main_missing_entries} summarize the main accuracy metrics across clean, noisy, ill-conditioned, and missing-entry regimes. Panel~\ref{fig:main_ablation_clean_error} summarizes the ablation behavior, while Panels~\ref{fig:main_ill_conditioned_kappaD}--\ref{fig:main_scaling_runtime} report core conditioning and runtime behavior. The figure shows the main empirical pattern developed in the tables below: the method behaves as expected on clean low-rank data, is competitive in noisy and ill-conditioned settings, and is not uniformly stronger in the missing-entry setting.

\subsection{Clean low-rank recovery}

Table~\ref{tab:full-clean-low-rank} reports the clean low-rank experiments. These results are a consistency check rather than a superiority test. Spectral methods are the appropriate benchmark for the unregularized full-data rank-constrained problem, and they recover clean low-rank targets to numerical precision when the target rank is sufficient. The proposed ridge and conditioned-core variants introduce a small regularization bias because their objective is penalized rather than purely unregularized. The relevant conclusion is that the proposed method behaves consistently with the rank-constrained formulation while also producing core-conditioning diagnostics.
\begin{table}[!htbp]
\centering
\tiny
\setlength{\tabcolsep}{2.5pt}
\renewcommand{\arraystretch}{0.92}
\caption{Clean low-rank recovery over the full multi-seed validation run. Spectral methods provide the reference for the unregularized full-data rank-constrained problem, while the proposed methods add regularization and core diagnostics.}
\label{tab:full-clean-low-rank}
\resizebox{\textwidth}{!}{%
\begin{tabular}{@{}llllllll@{}}
\toprule
Dataset & Method & RelFrob & CleanErr & $\kappa(D)$ & Time (s) & Sweeps & Storage \\
\midrule
clean\_m1000\_n700\_r020 & Truncated SVD & \(\num{5.77954145796e-15}\pm\num{1.49405244728e-15}\) & \(\num{5.77954145796e-15}\pm\num{1.49405244728e-15}\) & -- & \(\num{3.07914443999}\pm\num{0.838197510312}\) & \num{0} & \num{0.0491428571429} \\
clean\_m1000\_n700\_r020 & Randomized SVD & \num{4.28353387764e-15} & \num{4.28353387764e-15} & -- & \(\num{0.562094179995}\pm\num{0.820863364164}\) & \num{0} & \num{0.0491428571429} \\
clean\_m1000\_n700\_r020 & CUR & \(\num{2.29817044132e-14}\pm\num{1.5229652832e-14}\) & \(\num{2.29817044132e-14}\pm\num{1.5229652832e-14}\) & -- & \(\num{0.113237180014}\pm\num{0.0449485460838}\) & \num{0} & \num{0.0491428571429} \\
clean\_m1000\_n700\_r020 & Proposed ridge & \num{0.000829441245794} & \num{0.000829441245794} & \(\num{2.75423241164}\pm\num{9.41008164742e-15}\) & \(\num{47.01011198}\pm\num{7.22236474912}\) & \num{250} & \num{0.0491428571429} \\
clean\_m1000\_n700\_r020 & Proposed conditioned core & \num{0.000829441245794} & \num{0.000829441245794} & \(\num{2.75423241164}\pm\num{3.5943116778e-14}\) & \(\num{62.61270758}\pm\num{5.61234021669}\) & \num{250} & \num{0.0491428571429} \\
clean\_m1000\_n700\_r050 & Truncated SVD & \num{0.133513752929} & \num{0.133513752929} & -- & \(\num{2.87241712001}\pm\num{1.09064527572}\) & \num{0} & \num{0.0491428571429} \\
clean\_m1000\_n700\_r050 & Randomized SVD & \(\num{0.133534049312}\pm\num{4.15754291216e-06}\) & \(\num{0.133534049312}\pm\num{4.15754291216e-06}\) & -- & \(\num{0.128913299984}\pm\num{0.0501610821844}\) & \num{0} & \num{0.0491428571429} \\
clean\_m1000\_n700\_r050 & CUR & \(\num{0.226316257459}\pm\num{0.00309526732136}\) & \(\num{0.226316257459}\pm\num{0.00309526732136}\) & -- & \(\num{0.103168960009}\pm\num{0.027622699326}\) & \num{0} & \num{0.0491428571429} \\
clean\_m1000\_n700\_r050 & Proposed ridge & \(\num{0.133516283387}\pm\num{3.5653112107e-13}\) & \(\num{0.133516283387}\pm\num{3.5653112107e-13}\) & \(\num{2.7542260923}\pm\num{2.78364965017e-06}\) & \(\num{46.56359814}\pm\num{5.94247401265}\) & \num{250} & \num{0.0491428571429} \\
clean\_m1000\_n700\_r050 & Proposed conditioned core & \(\num{0.133516283387}\pm\num{3.56532775655e-13}\) & \(\num{0.133516283387}\pm\num{3.56532775655e-13}\) & \(\num{2.7542260923}\pm\num{2.78364963661e-06}\) & \(\num{64.51287868}\pm\num{8.08164136468}\) & \num{250} & \num{0.0491428571429} \\
clean\_m600\_n400\_r020 & Truncated SVD & \num{4.06489954994e-15} & \num{4.06489954994e-15} & -- & \(\num{1.23061820001}\pm\num{0.372977776568}\) & \num{0} & \num{0.085} \\
clean\_m600\_n400\_r020 & Randomized SVD & \(\num{4.60370504389e-15}\pm\num{1.88168209654e-15}\) & \(\num{4.60370504389e-15}\pm\num{1.88168209654e-15}\) & -- & \(\num{0.233494460001}\pm\num{0.165874430789}\) & \num{0} & \num{0.085} \\
clean\_m600\_n400\_r020 & CUR & \(\num{2.27323927048e-14}\pm\num{2.03874104797e-14}\) & \(\num{2.27323927048e-14}\pm\num{2.03874104797e-14}\) & -- & \(\num{0.150557559985}\pm\num{0.0652480888829}\) & \num{0} & \num{0.085} \\
clean\_m600\_n400\_r020 & Proposed ridge & \num{0.000829441245794} & \num{0.000829441245794} & \(\num{2.75423241164}\pm\num{1.66784914397e-14}\) & \(\num{34.76586442}\pm\num{4.36181843861}\) & \num{250} & \num{0.085} \\
clean\_m600\_n400\_r020 & Proposed conditioned core & \num{0.000829441245794} & \num{0.000829441245794} & \(\num{2.75423241164}\pm\num{2.26333062173e-14}\) & \(\num{50.16107182}\pm\num{10.4363250204}\) & \num{250} & \num{0.085} \\
clean\_m600\_n400\_r050 & Truncated SVD & \num{0.133513752929} & \num{0.133513752929} & -- & \(\num{0.895091960009}\pm\num{0.571249692848}\) & \num{0} & \num{0.085} \\
clean\_m600\_n400\_r050 & Randomized SVD & \(\num{0.133531166579}\pm\num{5.07350661013e-06}\) & \(\num{0.133531166579}\pm\num{5.07350661013e-06}\) & -- & \(\num{0.256258739985}\pm\num{0.345097898375}\) & \num{0} & \num{0.085} \\
clean\_m600\_n400\_r050 & CUR & \(\num{0.225120451683}\pm\num{0.00638219548649}\) & \(\num{0.225120451683}\pm\num{0.00638219548649}\) & -- & \(\num{0.104782740015}\pm\num{0.102070569988}\) & \num{0} & \num{0.085} \\
clean\_m600\_n400\_r050 & Proposed ridge & \(\num{0.133516283387}\pm\num{8.86835949979e-13}\) & \(\num{0.133516283387}\pm\num{8.86835949979e-13}\) & \(\num{2.75422675604}\pm\num{2.54115908019e-06}\) & \(\num{34.20898508}\pm\num{7.42690075204}\) & \num{250} & \num{0.085} \\
clean\_m600\_n400\_r050 & Proposed conditioned core & \(\num{0.133516283387}\pm\num{8.86837104121e-13}\) & \(\num{0.133516283387}\pm\num{8.86837104121e-13}\) & \(\num{2.75422675604}\pm\num{2.54115906958e-06}\) & \(\num{44.78370198}\pm\num{7.05139334749}\) & \num{250} & \num{0.085} \\
\bottomrule
\end{tabular}%
}
\end{table}

\subsection{Noisy low-rank recovery}

Table~\ref{tab:full-noisy-low-rank} reports noisy low-rank recovery across several signal-to-noise levels. The proposed variants track the spectral baselines closely in both relative Frobenius error and clean-target error. CUR is weaker in these runs, while ridge two-factor factorization is often much slower than the proposed tri-factorized updates. The main value of the proposed method in this regime is therefore not a large accuracy separation from SVD, but comparable recovery together with an explicit learned core and a reported \(\kappa(D)\). This supports the intended use of the method as a regularized and diagnostically transparent approximation rather than as a purely spectral replacement.
\begin{table}[!htbp]
\centering
\tiny
\setlength{\tabcolsep}{2.5pt}
\renewcommand{\arraystretch}{0.92}
\caption{Noisy low-rank recovery over the full multi-seed validation run. The proposed variants are interpreted against spectral and ridge baselines, with the additional diagnostic value of the learned core condition number.}
\label{tab:full-noisy-low-rank}
\resizebox{\textwidth}{!}{%
\begin{tabular}{@{}lllllll@{}}
\toprule
Dataset & Method & RelFrob & CleanErr & $\kappa(D)$ & Time (s) & Sweeps \\
\midrule
noisy\_m700\_n500\_r040\_snr10 & Truncated SVD & \(\num{0.219175264461}\pm\num{0.000185633573771}\) & \(\num{0.199626352836}\pm\num{7.05444695421e-06}\) & -- & \(\num{0.28990470001}\pm\num{0.159861764702}\) & \num{0} \\
noisy\_m700\_n500\_r040\_snr10 & Randomized SVD & \(\num{0.219175676109}\pm\num{0.000185301139317}\) & \(\num{0.199626822494}\pm\num{7.03566473498e-06}\) & -- & \(\num{0.037577179994}\pm\num{0.0145943266087}\) & \num{0} \\
noisy\_m700\_n500\_r040\_snr10 & CUR & \(\num{0.372648473151}\pm\num{0.00931737124623}\) & \(\num{0.36257149687}\pm\num{0.00976799343279}\) & -- & \(\num{0.0209179199999}\pm\num{0.00675204655138}\) & \num{0} \\
noisy\_m700\_n500\_r040\_snr10 & Ridge two-factor & \(\num{0.219176004159}\pm\num{0.000185554872466}\) & \(\num{0.1996249129}\pm\num{7.12854097947e-06}\) & -- & \(\num{225.15423784}\pm\num{93.6017078065}\) & \num{250} \\
noisy\_m700\_n500\_r040\_snr10 & Proposed ridge & \(\num{0.219175359495}\pm\num{0.000185633471528}\) & \(\num{0.199625409982}\pm\num{7.02444671031e-06}\) & \(\num{1.31092602912}\pm\num{0.00113114201897}\) & \(\num{9.57573239999}\pm\num{1.28222412603}\) & \num{250} \\
noisy\_m700\_n500\_r040\_snr10 & Proposed conditioned core & \(\num{0.219175359495}\pm\num{0.000185633471528}\) & \(\num{0.199625409982}\pm\num{7.02444671031e-06}\) & \(\num{1.31092602912}\pm\num{0.00113114201896}\) & \(\num{14.86994864}\pm\num{3.67015856402}\) & \num{250} \\
noisy\_m700\_n500\_r040\_snr2 & Truncated SVD & \(\num{0.466181050448}\pm\num{0.000249685146679}\) & \(\num{0.237740549764}\pm\num{0.000149860550343}\) & -- & \(\num{0.434382399998}\pm\num{0.44218182563}\) & \num{0} \\
noisy\_m700\_n500\_r040\_snr2 & Randomized SVD & \(\num{0.466183793893}\pm\num{0.000249043718507}\) & \(\num{0.237749744339}\pm\num{0.000151352776771}\) & -- & \(\num{0.059928799991}\pm\num{0.0547268336619}\) & \num{0} \\
noisy\_m700\_n500\_r040\_snr2 & CUR & \(\num{0.669777505593}\pm\num{0.00700115892996}\) & \(\num{0.573088981127}\pm\num{0.00928847238557}\) & -- & \(\num{0.024696060014}\pm\num{0.013872486734}\) & \num{0} \\
noisy\_m700\_n500\_r040\_snr2 & Ridge two-factor & \(\num{0.466181091533}\pm\num{0.000249684296004}\) & \(\num{0.237721122748}\pm\num{0.000149597296979}\) & -- & \(\num{427.576013}\pm\num{539.406482911}\) & \num{250} \\
noisy\_m700\_n500\_r040\_snr2 & Proposed ridge & \(\num{0.466181085654}\pm\num{0.000249685159597}\) & \(\num{0.237719395613}\pm\num{0.000149736015539}\) & \(\num{1.27990046646}\pm\num{0.00397742171298}\) & \(\num{13.76603768}\pm\num{10.7501246632}\) & \num{250} \\
noisy\_m700\_n500\_r040\_snr2 & Proposed conditioned core & \(\num{0.466181085654}\pm\num{0.000249685159597}\) & \(\num{0.237719395613}\pm\num{0.000149736015539}\) & \(\num{1.27990046646}\pm\num{0.00397742171298}\) & \(\num{19.76330152}\pm\num{11.2051383218}\) & \num{250} \\
noisy\_m700\_n500\_r040\_snr20 & Truncated SVD & \(\num{0.203445216361}\pm\num{8.60944558412e-05}\) & \(\num{0.198331062366}\pm\num{3.07434944379e-06}\) & -- & \(\num{0.247438439995}\pm\num{0.0781567478976}\) & \num{0} \\
noisy\_m700\_n500\_r040\_snr20 & Randomized SVD & \(\num{0.203445822849}\pm\num{8.65599886236e-05}\) & \(\num{0.198331724862}\pm\num{3.16219475714e-06}\) & -- & \(\num{0.0434810999897}\pm\num{0.0166053317083}\) & \num{0} \\
noisy\_m700\_n500\_r040\_snr20 & CUR & \(\num{0.332383566213}\pm\num{0.0111746489408}\) & \(\num{0.32947322801}\pm\num{0.0112703668715}\) & -- & \(\num{0.0310927999904}\pm\num{0.0216106297268}\) & \num{0} \\
noisy\_m700\_n500\_r040\_snr20 & Ridge two-factor & \(\num{0.203446127845}\pm\num{8.61331707365e-05}\) & \(\num{0.198331379497}\pm\num{3.08621747314e-06}\) & -- & \(\num{225.55357234}\pm\num{80.6780792327}\) & \num{250} \\
noisy\_m700\_n500\_r040\_snr20 & Proposed ridge & \(\num{0.203445319607}\pm\num{8.60944115265e-05}\) & \(\num{0.198330942948}\pm\num{3.08062628525e-06}\) & \(\num{1.31163860662}\pm\num{0.00076681138189}\) & \(\num{9.73649258}\pm\num{1.56323421006}\) & \num{250} \\
noisy\_m700\_n500\_r040\_snr20 & Proposed conditioned core & \(\num{0.203445319607}\pm\num{8.60944115265e-05}\) & \(\num{0.198330942948}\pm\num{3.08062628525e-06}\) & \(\num{1.31163860662}\pm\num{0.000766811381887}\) & \(\num{12.6354187}\pm\num{0.768515693111}\) & \num{250} \\
noisy\_m700\_n500\_r040\_snr5 & Truncated SVD & \(\num{0.271180341735}\pm\num{0.000122086618218}\) & \(\num{0.204671850872}\pm\num{5.45149796169e-05}\) & -- & \(\num{0.530411000014}\pm\num{0.507768608857}\) & \num{0} \\
noisy\_m700\_n500\_r040\_snr5 & Randomized SVD & \(\num{0.271180800842}\pm\num{0.00012241530454}\) & \(\num{0.204673005937}\pm\num{5.31382168511e-05}\) & -- & \(\num{0.0946928400081}\pm\num{0.0854305614186}\) & \num{0} \\
noisy\_m700\_n500\_r040\_snr5 & CUR & \(\num{0.455315324865}\pm\num{0.00929162905237}\) & \(\num{0.423849082229}\pm\num{0.0104068948306}\) & -- & \(\num{0.0379937999882}\pm\num{0.0172996446673}\) & \num{0} \\
noisy\_m700\_n500\_r040\_snr5 & Ridge two-factor & \(\num{0.2711807128}\pm\num{0.000122117306529}\) & \(\num{0.204665423266}\pm\num{5.46899804463e-05}\) & -- & \(\num{381.40875322}\pm\num{406.774848754}\) & \num{250} \\
noisy\_m700\_n500\_r040\_snr5 & Proposed ridge & \(\num{0.271180416139}\pm\num{0.000122086589977}\) & \(\num{0.204668188277}\pm\num{5.4432796898e-05}\) & \(\num{1.30708514811}\pm\num{0.00192883465477}\) & \(\num{10.49150328}\pm\num{2.12452692575}\) & \num{250} \\
noisy\_m700\_n500\_r040\_snr5 & Proposed conditioned core & \(\num{0.271180416139}\pm\num{0.000122086589977}\) & \(\num{0.204668188277}\pm\num{5.4432796898e-05}\) & \(\num{1.30708514811}\pm\num{0.00192883465471}\) & \(\num{14.69633166}\pm\num{2.09361495829}\) & \num{250} \\
noisy\_m700\_n500\_r040\_snr50 & Truncated SVD & \(\num{0.198808898323}\pm\num{3.69253235365e-05}\) & \(\num{0.197971875119}\pm\num{6.68902825081e-07}\) & -- & \(\num{0.188319499977}\pm\num{0.0182366791831}\) & \num{0} \\
noisy\_m700\_n500\_r040\_snr50 & Randomized SVD & \(\num{0.198809550036}\pm\num{3.69138918087e-05}\) & \(\num{0.197972519957}\pm\num{1.62248156618e-06}\) & -- & \(\num{0.0346533200005}\pm\num{0.011250907481}\) & \num{0} \\
noisy\_m700\_n500\_r040\_snr50 & CUR & \(\num{0.323476711626}\pm\num{0.0106566106338}\) & \(\num{0.32298335917}\pm\num{0.0106713303606}\) & -- & \(\num{0.0190332799917}\pm\num{0.00735464106832}\) & \num{0} \\
noisy\_m700\_n500\_r040\_snr50 & Ridge two-factor & \(\num{0.198809868692}\pm\num{3.68895394665e-05}\) & \(\num{0.197972832107}\pm\num{7.43185800473e-07}\) & -- & \(\num{190.00186838}\pm\num{5.75187973259}\) & \num{250} \\
noisy\_m700\_n500\_r040\_snr50 & Proposed ridge & \(\num{0.198809004218}\pm\num{3.69253000461e-05}\) & \(\num{0.19797195627}\pm\num{6.59710802853e-07}\) & \(\num{1.31242109606}\pm\num{0.000121421729137}\) & \(\num{9.08496448001}\pm\num{0.4194501952}\) & \num{250} \\
noisy\_m700\_n500\_r040\_snr50 & Proposed conditioned core & \(\num{0.198809004218}\pm\num{3.6925300046e-05}\) & \(\num{0.19797195627}\pm\num{6.59710802839e-07}\) & \(\num{1.31242109606}\pm\num{0.000121421729183}\) & \(\num{13.39911754}\pm\num{1.22172689562}\) & \num{250} \\
\bottomrule
\end{tabular}%
}
\end{table}

\subsection{Ill-conditioned spectra}

Table~\ref{tab:full-ill-conditioned} evaluates matrices with increasingly severe spectral profiles. This is one of the most relevant tests for the proposed method because core conditioning is part of the model design. The results show that the proposed variants remain competitive in clean-target error. In the most ill-conditioned cases, the unconstrained ridge core may become singular, reported as \(\kappa(D)=\infty\). The conditioned-core variant prevents this by enforcing a finite admissible core condition number. This supports the central claim that conditioning of the latent core should be treated as an explicit modeling and diagnostic quantity rather than as an automatic consequence of the factorization.
\begin{table}[!htbp]
\centering
\tiny
\setlength{\tabcolsep}{2.5pt}
\renewcommand{\arraystretch}{0.92}
\caption{Ill-conditioned spectra over the full multi-seed validation run. These experiments test how the approximation and learned core behave as the spectral condition profile becomes more severe.}
\label{tab:full-ill-conditioned}
\resizebox{\textwidth}{!}{%
\begin{tabular}{@{}lllllll@{}}
\toprule
Dataset & Method & RelFrob & CleanErr & $\kappa(D)$ & Time (s) & Sweeps \\
\midrule
illcond\_m600\_eta10 & Truncated SVD & \(\num{0.019173770046}\pm\num{4.44911912618e-06}\) & \(\num{0.0058773863218}\pm\num{1.61130875336e-05}\) & -- & \(\num{0.380544540007}\pm\num{0.0806139663784}\) & \num{0} \\
illcond\_m600\_eta10 & Randomized SVD & \(\num{0.0192089790984}\pm\num{6.09461557135e-06}\) & \(\num{0.00581117735225}\pm\num{1.69590319399e-05}\) & -- & \(\num{0.0491424800013}\pm\num{0.0193611962199}\) & \num{0} \\
illcond\_m600\_eta10 & CUR & \(\num{0.030037471001}\pm\num{0.000780406554825}\) & \(\num{0.0230306723401}\pm\num{0.00101393950012}\) & -- & \(\num{0.0298529999913}\pm\num{0.0144991093649}\) & \num{0} \\
illcond\_m600\_eta10 & Proposed ridge & \(\num{0.0196042816645}\pm\num{6.03108454177e-06}\) & \(\num{0.0045305760295}\pm\num{2.13132811469e-05}\) & $\infty$ & \(\num{10.91594598}\pm\num{1.94526039258}\) & \num{250} \\
illcond\_m600\_eta10 & Proposed conditioned core & \(\num{0.0196042816645}\pm\num{6.03108454177e-06}\) & \(\num{0.0045305760295}\pm\num{2.13132811469e-05}\) & \(\num{999.994932707}\pm\num{0.000543364036203}\) & \(\num{14.26453852}\pm\num{1.41502787356}\) & \num{250} \\
illcond\_m600\_eta2 & Truncated SVD & \(\num{0.15354947238}\pm\num{2.07803036763e-05}\) & \(\num{0.152462816769}\pm\num{3.02680571692e-07}\) & -- & \(\num{0.26993503999}\pm\num{0.0937684081827}\) & \num{0} \\
illcond\_m600\_eta2 & Randomized SVD & \(\num{0.153549494774}\pm\num{2.07861836388e-05}\) & \(\num{0.152462840087}\pm\num{3.03280764989e-07}\) & -- & \(\num{0.0336733999895}\pm\num{0.00515713024679}\) & \num{0} \\
illcond\_m600\_eta2 & CUR & \(\num{0.259143577494}\pm\num{0.0118603099903}\) & \(\num{0.258525960656}\pm\num{0.0119134128859}\) & -- & \(\num{0.0196445799898}\pm\num{0.00493826764995}\) & \num{0} \\
illcond\_m600\_eta2 & Proposed ridge & \(\num{0.153549661567}\pm\num{2.07802785606e-05}\) & \(\num{0.152462941493}\pm\num{3.17380510752e-07}\) & \(\num{1.50008050971}\pm\num{0.000324022808009}\) & \(\num{9.69833102002}\pm\num{1.02108891107}\) & \num{250} \\
illcond\_m600\_eta2 & Proposed conditioned core & \(\num{0.153549661567}\pm\num{2.07802785607e-05}\) & \(\num{0.152462941493}\pm\num{3.17380510754e-07}\) & \(\num{1.50008050971}\pm\num{0.000324022808012}\) & \(\num{14.0356979}\pm\num{0.528347808841}\) & \num{250} \\
illcond\_m600\_eta4 & Truncated SVD & \(\num{0.0302827312628}\pm\num{1.94720335639e-05}\) & \(\num{0.0238528121577}\pm\num{4.07045635209e-06}\) & -- & \(\num{0.29543326}\pm\num{0.0700125652117}\) & \num{0} \\
illcond\_m600\_eta4 & Randomized SVD & \(\num{0.0302827312632}\pm\num{1.94720342041e-05}\) & \(\num{0.023852812264}\pm\num{4.07053057768e-06}\) & -- & \(\num{0.0329334599781}\pm\num{0.00581946746681}\) & \num{0} \\
illcond\_m600\_eta4 & CUR & \(\num{0.0682702655691}\pm\num{0.00490170617176}\) & \(\num{0.0656538905918}\pm\num{0.00508266947126}\) & -- & \(\num{0.0209641799913}\pm\num{0.00784609882976}\) & \num{0} \\
illcond\_m600\_eta4 & Proposed ridge & \(\num{0.0302874380938}\pm\num{1.94690670716e-05}\) & \(\num{0.0238557541045}\pm\num{4.20853605852e-06}\) & \(\num{3.29892047999}\pm\num{0.00174414726255}\) & \(\num{9.91710706}\pm\num{0.789830621171}\) & \num{250} \\
illcond\_m600\_eta4 & Proposed conditioned core & \(\num{0.0302874380938}\pm\num{1.94690670716e-05}\) & \(\num{0.0238557541045}\pm\num{4.20853605852e-06}\) & \(\num{3.29892047999}\pm\num{0.00174414726261}\) & \(\num{14.7745939}\pm\num{1.25571154662}\) & \num{250} \\
illcond\_m600\_eta6 & Truncated SVD & \(\num{0.0196573873989}\pm\num{9.27693798717e-06}\) & \(\num{0.00625539759459}\pm\num{1.49641912059e-05}\) & -- & \(\num{0.31183784001}\pm\num{0.0929870361825}\) & \num{0} \\
illcond\_m600\_eta6 & Randomized SVD & \(\num{0.0196574068595}\pm\num{9.28123300706e-06}\) & \(\num{0.00625547898518}\pm\num{1.49271303378e-05}\) & -- & \(\num{0.0359966399963}\pm\num{0.0142277511983}\) & \num{0} \\
illcond\_m600\_eta6 & CUR & \(\num{0.0410636235111}\pm\num{0.00243260668035}\) & \(\num{0.0364405693455}\pm\num{0.00270909320425}\) & -- & \(\num{0.0247314200155}\pm\num{0.0143844164282}\) & \num{0} \\
illcond\_m600\_eta6 & Proposed ridge & \(\num{0.0196805363617}\pm\num{9.26333947699e-06}\) & \(\num{0.00624450595814}\pm\num{1.36134915481e-05}\) & \(\num{6.09117325074}\pm\num{0.0124612937017}\) & \(\num{10.11594798}\pm\num{1.22218051639}\) & \num{250} \\
illcond\_m600\_eta6 & Proposed conditioned core & \(\num{0.0196805363617}\pm\num{9.26333947699e-06}\) & \(\num{0.00624450595814}\pm\num{1.36134915481e-05}\) & \(\num{6.09117325074}\pm\num{0.0124612937017}\) & \(\num{14.72812082}\pm\num{1.19622008906}\) & \num{250} \\
illcond\_m600\_eta8 & Truncated SVD & \(\num{0.019284949604}\pm\num{4.02470795351e-06}\) & \(\num{0.00558468407054}\pm\num{2.78807044054e-05}\) & -- & \(\num{0.294373699999}\pm\num{0.0986515926602}\) & \num{0} \\
illcond\_m600\_eta8 & Randomized SVD & \(\num{0.0193022237038}\pm\num{5.22921190322e-06}\) & \(\num{0.00559351605144}\pm\num{1.76231421114e-05}\) & -- & \(\num{0.0313405200141}\pm\num{0.00593956859811}\) & \num{0} \\
illcond\_m600\_eta8 & CUR & \(\num{0.0327241037673}\pm\num{0.00118958171958}\) & \(\num{0.0265556780776}\pm\num{0.00146084481854}\) & -- & \(\num{0.0178165400052}\pm\num{0.0028663216777}\) & \num{0} \\
illcond\_m600\_eta8 & Proposed ridge & \(\num{0.0194620263113}\pm\num{3.70011608925e-06}\) & \(\num{0.00499063061004}\pm\num{2.18155586019e-05}\) & $\infty$ & \(\num{10.02473846}\pm\num{0.691713734949}\) & \num{250} \\
illcond\_m600\_eta8 & Proposed conditioned core & \(\num{0.0194620263113}\pm\num{3.70011608925e-06}\) & \(\num{0.00499063061004}\pm\num{2.18155586019e-05}\) & \(\num{999.995701044}\pm\num{0.000783258584031}\) & \(\num{15.1557252}\pm\num{1.39232116515}\) & \num{250} \\
\bottomrule
\end{tabular}%
}
\end{table}

\subsection{Missing-entry recovery}

Table~\ref{tab:full-missing-entries} reports the missing-entry experiments. The relevant quantities in this setting are validation and test RMSE, not training RMSE alone. The proposed weighted tri-factorization fits the observed entries closely, but the held-out RMSE shows that it is not uniformly better than the matrix-completion baselines. In particular, soft-impute and ridge two-factor factorization are stronger in several observed-rate regimes. At the highest observed rate, the proposed variants become competitive, but the overall evidence does not justify a broad missing-entry superiority claim. We therefore treat the weighted formulation as a useful extension of the framework, not as the main empirical strength of the present paper.
\begin{table}[!htbp]
\centering
\tiny
\setlength{\tabcolsep}{2.5pt}
\renewcommand{\arraystretch}{0.92}
\caption{Missing-entry recovery over the full multi-seed validation run. Held-out validation and test RMSE are the relevant generalization metrics for this setting.}
\label{tab:full-missing-entries}
\resizebox{\textwidth}{!}{%
\begin{tabular}{@{}lllllllll@{}}
\toprule
Dataset & Method & Train RMSE & Val RMSE & Test RMSE & RelFrob & CleanErr & $\kappa(D)$ & Time (s) \\
\midrule
missing\_m450\_n320\_obs0.05 & Mean-imputed SVD & \(\num{0.00283805985588}\pm\num{5.71677636682e-05}\) & \(\num{0.00371952783244}\pm\num{8.74255742364e-05}\) & \(\num{0.00384684818388}\pm\num{0.000122301627662}\) & \(\num{0.983535166709}\pm\num{0.000428645511983}\) & \(\num{0.983525636048}\pm\num{0.000432233752167}\) & -- & \(\num{0.0796426399844}\pm\num{0.0286122647205}\) \\
missing\_m450\_n320\_obs0.05 & Soft-impute & \(\num{2.29363973988e-05}\pm\num{2.83351001873e-07}\) & \(\num{0.00378355142199}\pm\num{0.000124221205419}\) & \(\num{0.00392394841355}\pm\num{0.000136021186278}\) & \(\num{0.987030277107}\pm\num{0.00408144819037}\) & \(\num{0.987054252093}\pm\num{0.00406547130098}\) & -- & \(\num{14.04470786}\pm\num{6.02059071212}\) \\
missing\_m450\_n320\_obs0.05 & Ridge two-factor & \(\num{1.83792093358e-05}\pm\num{9.21060743378e-08}\) & \(\num{0.00322193839497}\pm\num{7.93153028884e-05}\) & \(\num{0.00324902713638}\pm\num{9.85555628434e-05}\) & \(\num{0.830345774459}\pm\num{0.00932939119874}\) & \(\num{0.830266760675}\pm\num{0.00933804472101}\) & -- & \(\num{26.40327246}\pm\num{9.81029022408}\) \\
missing\_m450\_n320\_obs0.05 & Proposed ridge & \(\num{3.91193957147e-05}\pm\num{3.43650481604e-07}\) & \(\num{0.00339918070124}\pm\num{5.05277988082e-05}\) & \(\num{0.00346553059581}\pm\num{0.00015284459154}\) & \(\num{0.89392319403}\pm\num{0.0177680346184}\) & \(\num{0.893915168362}\pm\num{0.0177840393505}\) & $\infty$ & \(\num{1575.5152546}\pm\num{426.598445121}\) \\
missing\_m450\_n320\_obs0.05 & Proposed conditioned core & \(\num{3.91193957147e-05}\pm\num{3.43650481621e-07}\) & \(\num{0.00339918070124}\pm\num{5.05277988073e-05}\) & \(\num{0.00346553059581}\pm\num{0.00015284459154}\) & \(\num{0.89392319403}\pm\num{0.0177680346182}\) & \(\num{0.893915168362}\pm\num{0.0177840393503}\) & \(\num{11659.4619352}\pm\num{310.557788534}\) & \(\num{1570.05818926}\pm\num{484.779171799}\) \\
missing\_m450\_n320\_obs0.10 & Mean-imputed SVD & \(\num{0.00278647704135}\pm\num{2.01154606924e-05}\) & \(\num{0.00372648009991}\pm\num{4.5779242404e-05}\) & \(\num{0.0036979554031}\pm\num{5.36278233824e-05}\) & \(\num{0.95990055027}\pm\num{0.000885207296784}\) & \(\num{0.959888972527}\pm\num{0.000885095732612}\) & -- & \(\num{0.0616076399805}\pm\num{0.0237013261578}\) \\
missing\_m450\_n320\_obs0.10 & Soft-impute & \(\num{8.99979420435e-05}\pm\num{3.47331071755e-06}\) & \(\num{0.00369287502548}\pm\num{3.98630385708e-05}\) & \(\num{0.0036836191109}\pm\num{4.35895681416e-05}\) & \(\num{0.938472499211}\pm\num{0.00727638005372}\) & \(\num{0.938481070181}\pm\num{0.00725030387409}\) & -- & \(\num{15.81512302}\pm\num{4.99928190179}\) \\
missing\_m450\_n320\_obs0.10 & Ridge two-factor & \(\num{2.2237408424e-05}\pm\num{4.84686825634e-07}\) & \(\num{0.00337136507348}\pm\num{0.000102566124343}\) & \(\num{0.00334028339488}\pm\num{8.20726186984e-05}\) & \(\num{0.849561529367}\pm\num{0.0154454314043}\) & \(\num{0.849542130526}\pm\num{0.015512015818}\) & -- & \(\num{27.50090306}\pm\num{14.4655364219}\) \\
missing\_m450\_n320\_obs0.10 & Proposed ridge & \(\num{3.89324837037e-05}\pm\num{7.70089402926e-07}\) & \(\num{0.00372419549759}\pm\num{7.93843396407e-05}\) & \(\num{0.00372061363063}\pm\num{6.80385780429e-05}\) & \(\num{0.947604257382}\pm\num{0.0222275607137}\) & \(\num{0.947577244405}\pm\num{0.0222457157425}\) & \(\num{2.29398389133}\pm\num{0.137943880654}\) & \(\num{3130.87412818}\pm\num{1441.38204241}\) \\
missing\_m450\_n320\_obs0.10 & Proposed conditioned core & \(\num{3.89324837037e-05}\pm\num{7.70089402926e-07}\) & \(\num{0.00372419549759}\pm\num{7.93843396407e-05}\) & \(\num{0.00372061363063}\pm\num{6.80385780429e-05}\) & \(\num{0.947604257382}\pm\num{0.0222275607137}\) & \(\num{0.947577244405}\pm\num{0.0222457157425}\) & \(\num{2.29398389133}\pm\num{0.137943880654}\) & \(\num{2900.34800302}\pm\num{1536.34025463}\) \\
missing\_m450\_n320\_obs0.20 & Mean-imputed SVD & \(\num{0.00267777714947}\pm\num{1.23153677721e-05}\) & \(\num{0.00356840358608}\pm\num{5.36493163072e-05}\) & \(\num{0.00357204248602}\pm\num{6.9948334136e-05}\) & \(\num{0.908648018378}\pm\num{0.00120515985751}\) & \(\num{0.908611568355}\pm\num{0.00121034064716}\) & -- & \(\num{0.0438133600052}\pm\num{0.0295431833373}\) \\
missing\_m450\_n320\_obs0.20 & Soft-impute & \(\num{0.000477601532001}\pm\num{6.08792197661e-06}\) & \(\num{0.00282912705197}\pm\num{9.27206061242e-05}\) & \(\num{0.00279777627876}\pm\num{5.96694790836e-05}\) & \(\num{0.689417041741}\pm\num{0.0163622861872}\) & \(\num{0.689260176482}\pm\num{0.0163606687571}\) & -- & \(\num{12.17088216}\pm\num{6.57402361297}\) \\
missing\_m450\_n320\_obs0.20 & Ridge two-factor & \(\num{0.000364248857858}\pm\num{1.33718970088e-05}\) & \(\num{0.00394106858028}\pm\num{0.000160960922902}\) & \(\num{0.0039606812561}\pm\num{0.000198535412116}\) & \(\num{0.96003605302}\pm\num{0.0356176174603}\) & \(\num{0.960048425382}\pm\num{0.0356277754882}\) & -- & \(\num{20.05118836}\pm\num{11.2808254421}\) \\
missing\_m450\_n320\_obs0.20 & Proposed ridge & \(\num{0.000385550914945}\pm\num{1.3251511529e-05}\) & \(\num{0.00369880449033}\pm\num{0.000190814107874}\) & \(\num{0.00369538167485}\pm\num{0.000168506158513}\) & \(\num{0.895718131638}\pm\num{0.0461379783278}\) & \(\num{0.89571246363}\pm\num{0.0461661971969}\) & \(\num{2.50807056144}\pm\num{0.0515441027167}\) & \(\num{5307.0418633}\pm\num{2858.98465201}\) \\
missing\_m450\_n320\_obs0.20 & Proposed conditioned core & \(\num{0.000385550914945}\pm\num{1.3251511529e-05}\) & \(\num{0.00369880449033}\pm\num{0.000190814107874}\) & \(\num{0.00369538167485}\pm\num{0.000168506158513}\) & \(\num{0.895718131638}\pm\num{0.0461379783278}\) & \(\num{0.89571246363}\pm\num{0.0461661971969}\) & \(\num{2.50807056144}\pm\num{0.0515441027168}\) & \(\num{7248.3970922}\pm\num{6361.87535876}\) \\
missing\_m450\_n320\_obs0.40 & Mean-imputed SVD & \(\num{0.0023876706295}\pm\num{6.93119690295e-06}\) & \(\num{0.00323887706246}\pm\num{7.91873334114e-05}\) & \(\num{0.00324822593762}\pm\num{3.11345579369e-05}\) & \(\num{0.796066205518}\pm\num{0.0032130562639}\) & \(\num{0.79595744167}\pm\num{0.00320902881981}\) & -- & \(\num{0.065844899998}\pm\num{0.0428578191877}\) \\
missing\_m450\_n320\_obs0.40 & Soft-impute & \(\num{0.000532800791354}\pm\num{5.72571896631e-06}\) & \(\num{0.000905706898655}\pm\num{3.57963524449e-05}\) & \(\num{0.000899542217066}\pm\num{2.73543742927e-05}\) & \(\num{0.214937485514}\pm\num{0.00651894589192}\) & \(\num{0.213839430125}\pm\num{0.00655857208609}\) & -- & \(\num{15.8576427}\pm\num{12.0946967368}\) \\
missing\_m450\_n320\_obs0.40 & Ridge two-factor & \(\num{0.000528013136349}\pm\num{2.66235115783e-06}\) & \(\num{0.000893101666243}\pm\num{1.19890487913e-05}\) & \(\num{0.00088670757955}\pm\num{3.25484157321e-06}\) & \(\num{0.211780931916}\pm\num{0.00117572686859}\) & \(\num{0.21068916202}\pm\num{0.00120001021163}\) & -- & \(\num{31.33951936}\pm\num{26.3512841182}\) \\
missing\_m450\_n320\_obs0.40 & Proposed ridge & \(\num{0.000528034378882}\pm\num{2.46581660586e-06}\) & \(\num{0.000890344575481}\pm\num{1.27903722059e-05}\) & \(\num{0.000882873002104}\pm\num{5.44383682697e-06}\) & \(\num{0.210931259865}\pm\num{0.000930862086007}\) & \(\num{0.209829219785}\pm\num{0.00101390763156}\) & \(\num{2.49147648547}\pm\num{0.0416616389692}\) & \(\num{11574.8551574}\pm\num{5482.15077794}\) \\
missing\_m450\_n320\_obs0.40 & Proposed conditioned core & \(\num{0.000528034378882}\pm\num{2.46581660586e-06}\) & \(\num{0.000890344575481}\pm\num{1.27903722059e-05}\) & \(\num{0.000882873002104}\pm\num{5.44383682697e-06}\) & \(\num{0.210931259865}\pm\num{0.000930862086007}\) & \(\num{0.209829219785}\pm\num{0.00101390763156}\) & \(\num{2.49147648548}\pm\num{0.0416616389692}\) & \(\num{10942.4028913}\pm\num{1813.0975385}\) \\
\bottomrule
\end{tabular}%
}
\end{table}

\subsection{Convergence and ablation diagnostics}

\begin{figure}[!htbp]
\centering
\begin{subfigure}[t]{0.24\textwidth}
\includegraphics[width=\linewidth]{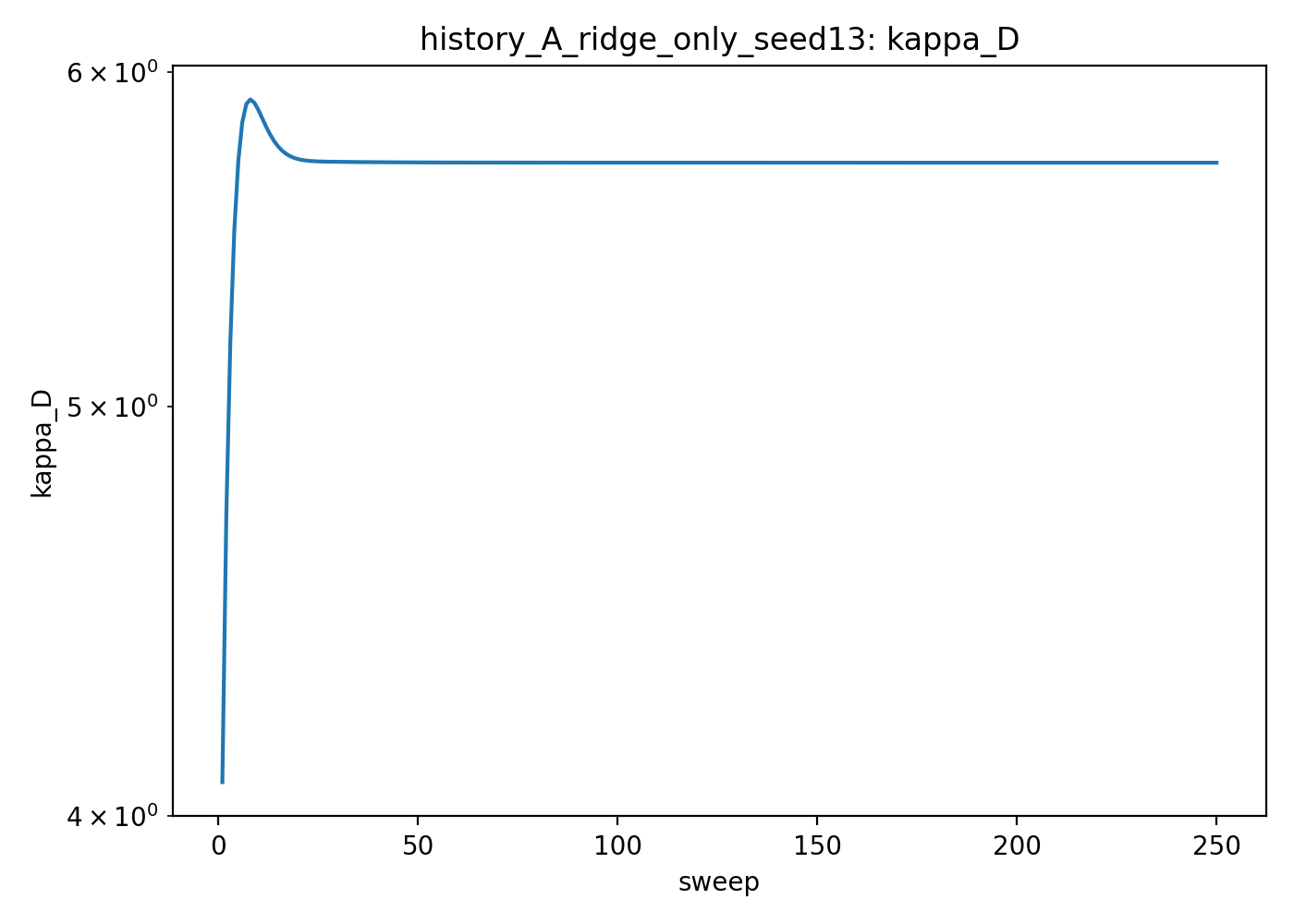}
\caption{A: \(\kappa(D)\).}
\label{fig:main_A_kappa_seed13}
\end{subfigure}\hfill
\begin{subfigure}[t]{0.24\textwidth}
\includegraphics[width=\linewidth]{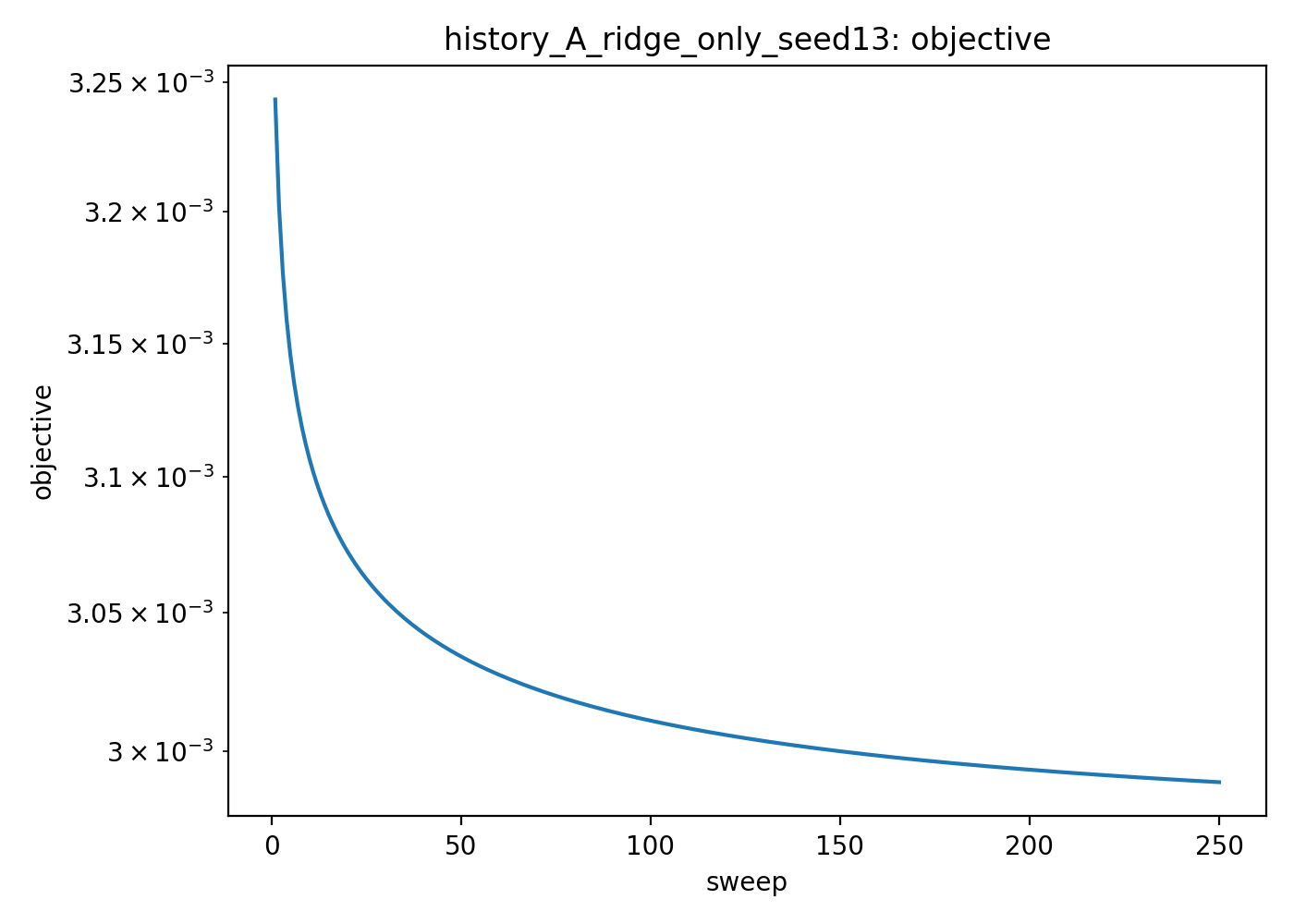}
\caption{A: objective.}
\label{fig:main_A_objective_seed13}
\end{subfigure}\hfill
\begin{subfigure}[t]{0.24\textwidth}
\includegraphics[width=\linewidth]{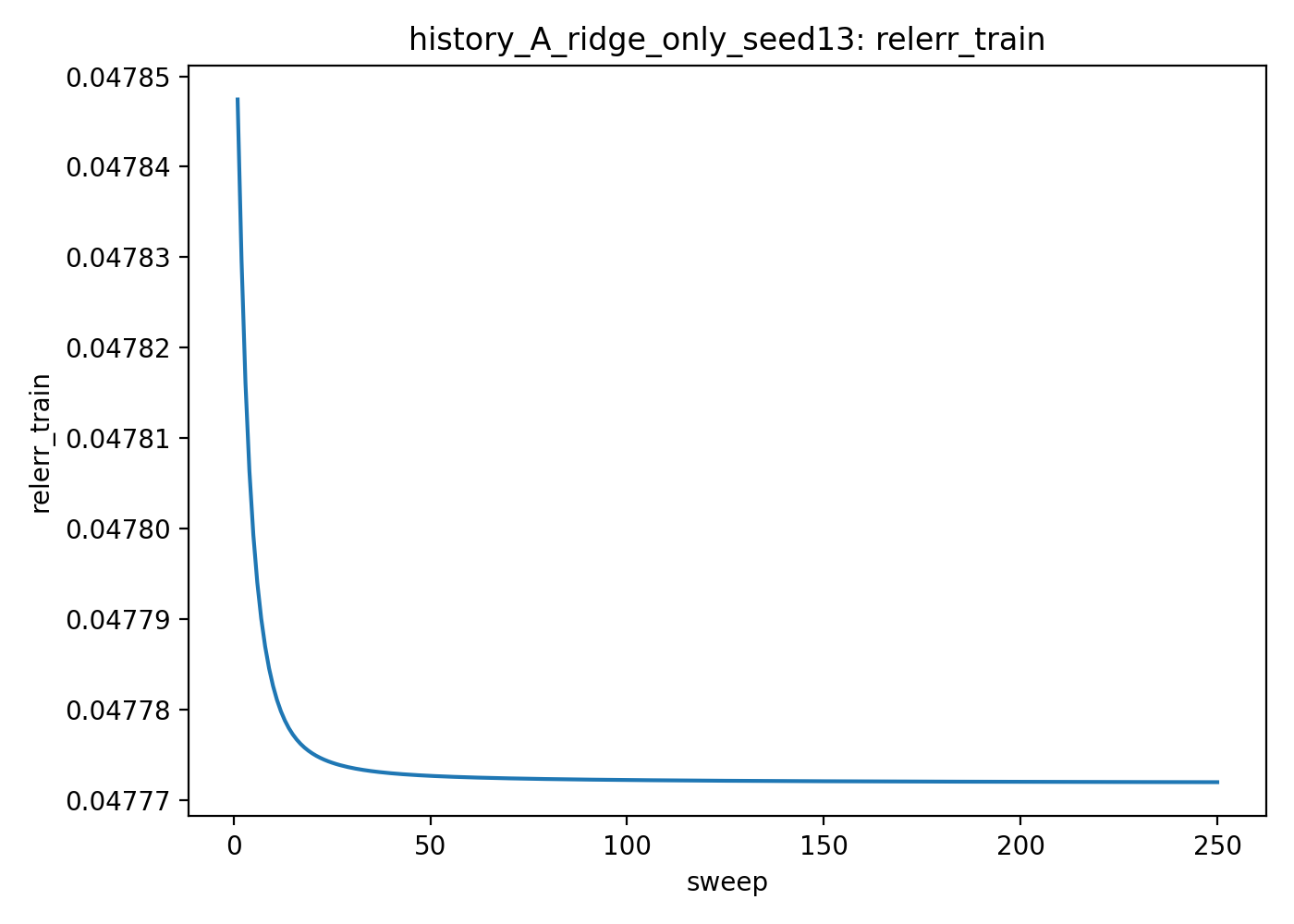}
\caption{A: training residual.}
\label{fig:main_A_residual_seed13}
\end{subfigure}\hfill
\begin{subfigure}[t]{0.24\textwidth}
\includegraphics[width=\linewidth]{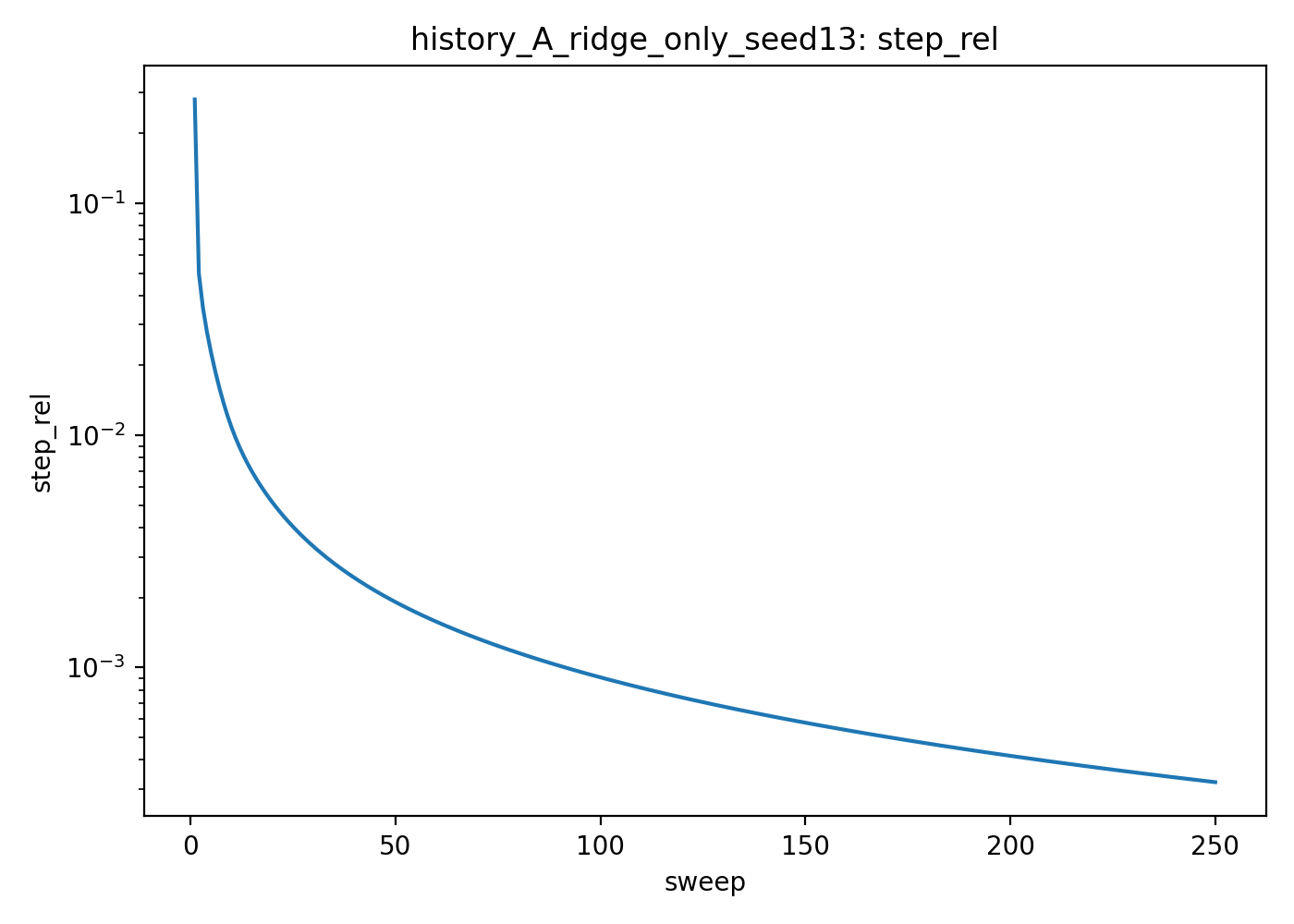}
\caption{A: relative step.}
\label{fig:main_A_step_seed13}
\end{subfigure}

\medskip

\begin{subfigure}[t]{0.24\textwidth}
\includegraphics[width=\linewidth]{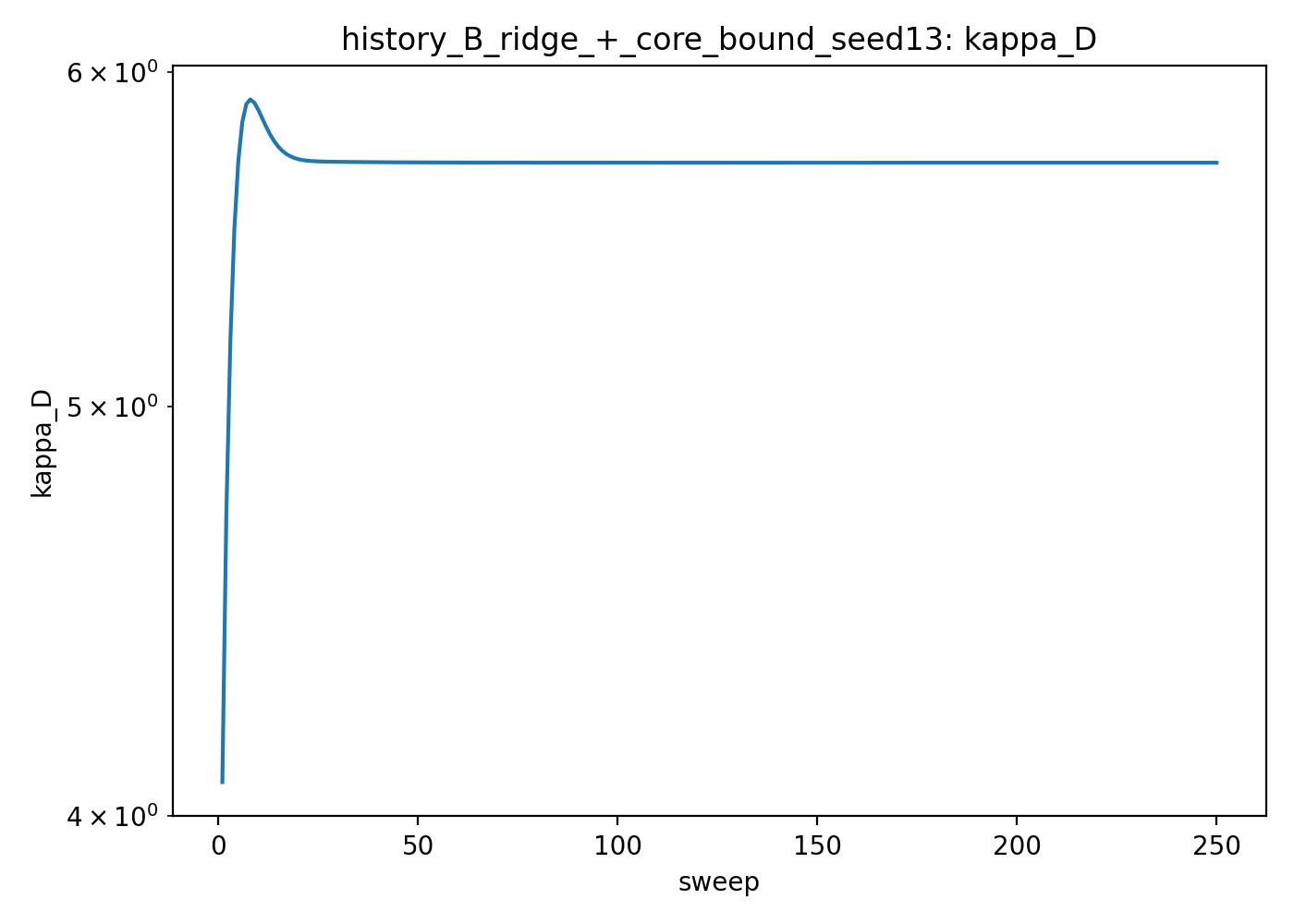}
\caption{B: \(\kappa(D)\).}
\label{fig:main_B_kappa_seed13}
\end{subfigure}\hfill
\begin{subfigure}[t]{0.24\textwidth}
\includegraphics[width=\linewidth]{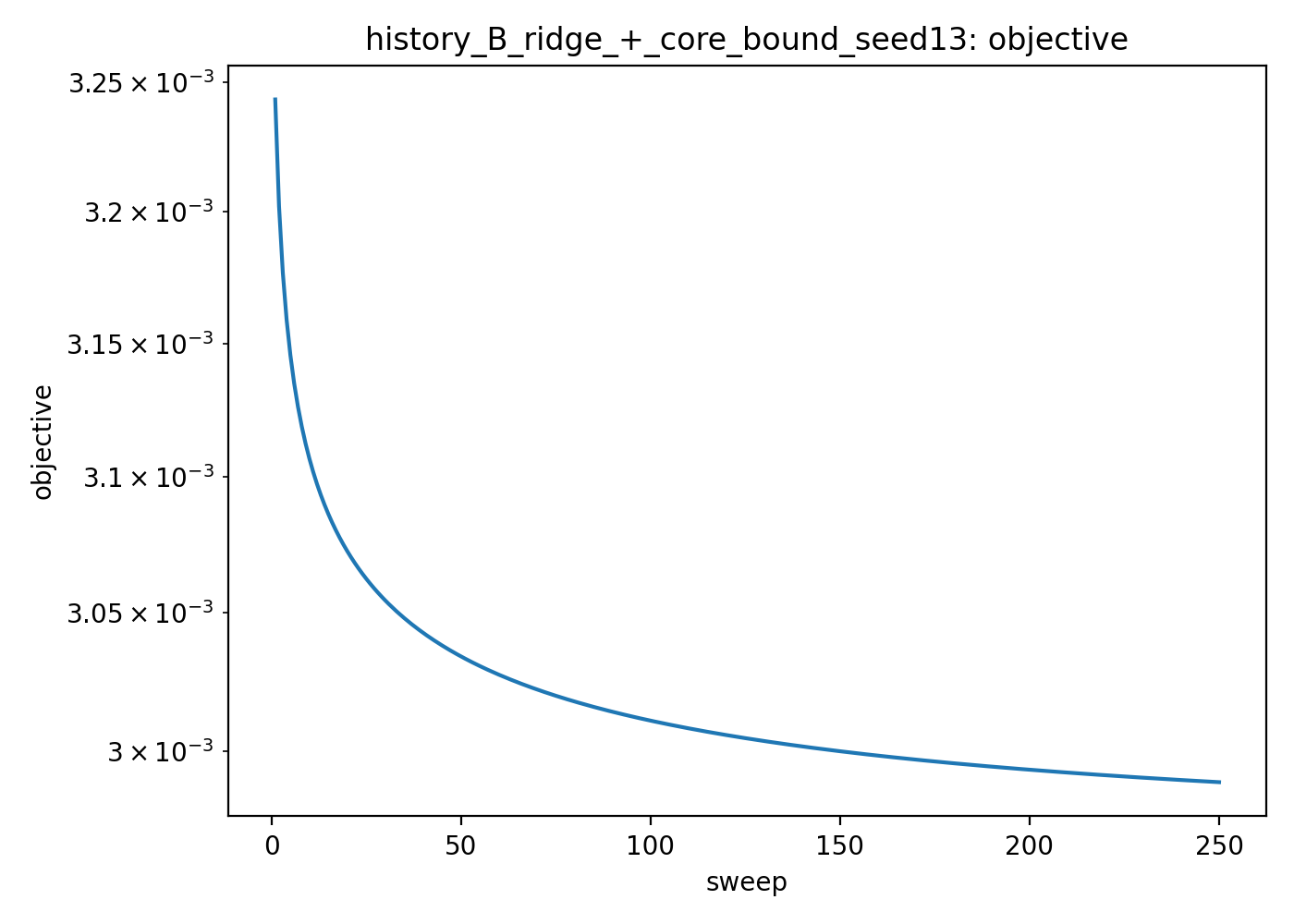}
\caption{B: objective.}
\label{fig:main_B_objective_seed13}
\end{subfigure}\hfill
\begin{subfigure}[t]{0.24\textwidth}
\includegraphics[width=\linewidth]{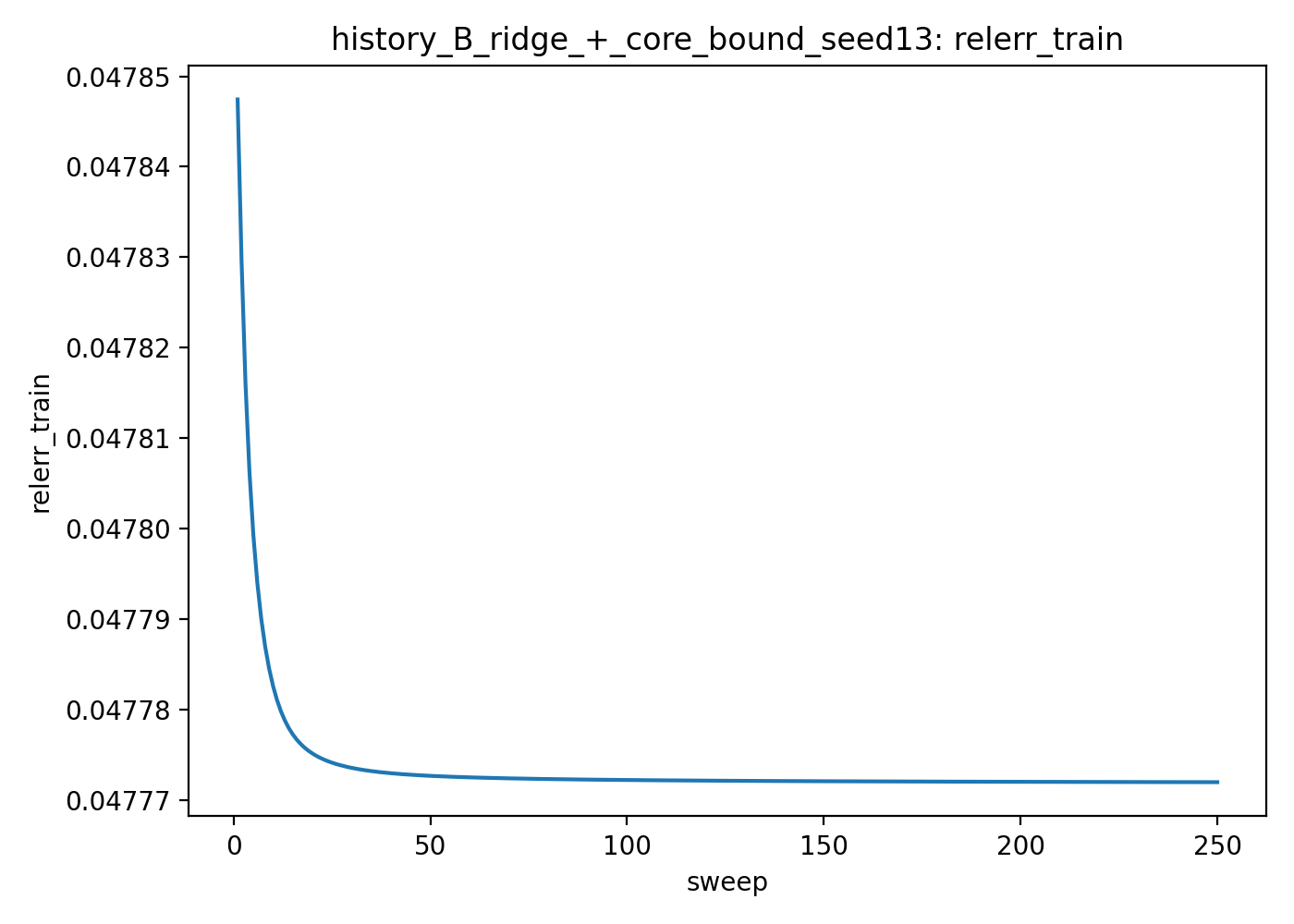}
\caption{B: training residual.}
\label{fig:main_B_residual_seed13}
\end{subfigure}\hfill
\begin{subfigure}[t]{0.24\textwidth}
\includegraphics[width=\linewidth]{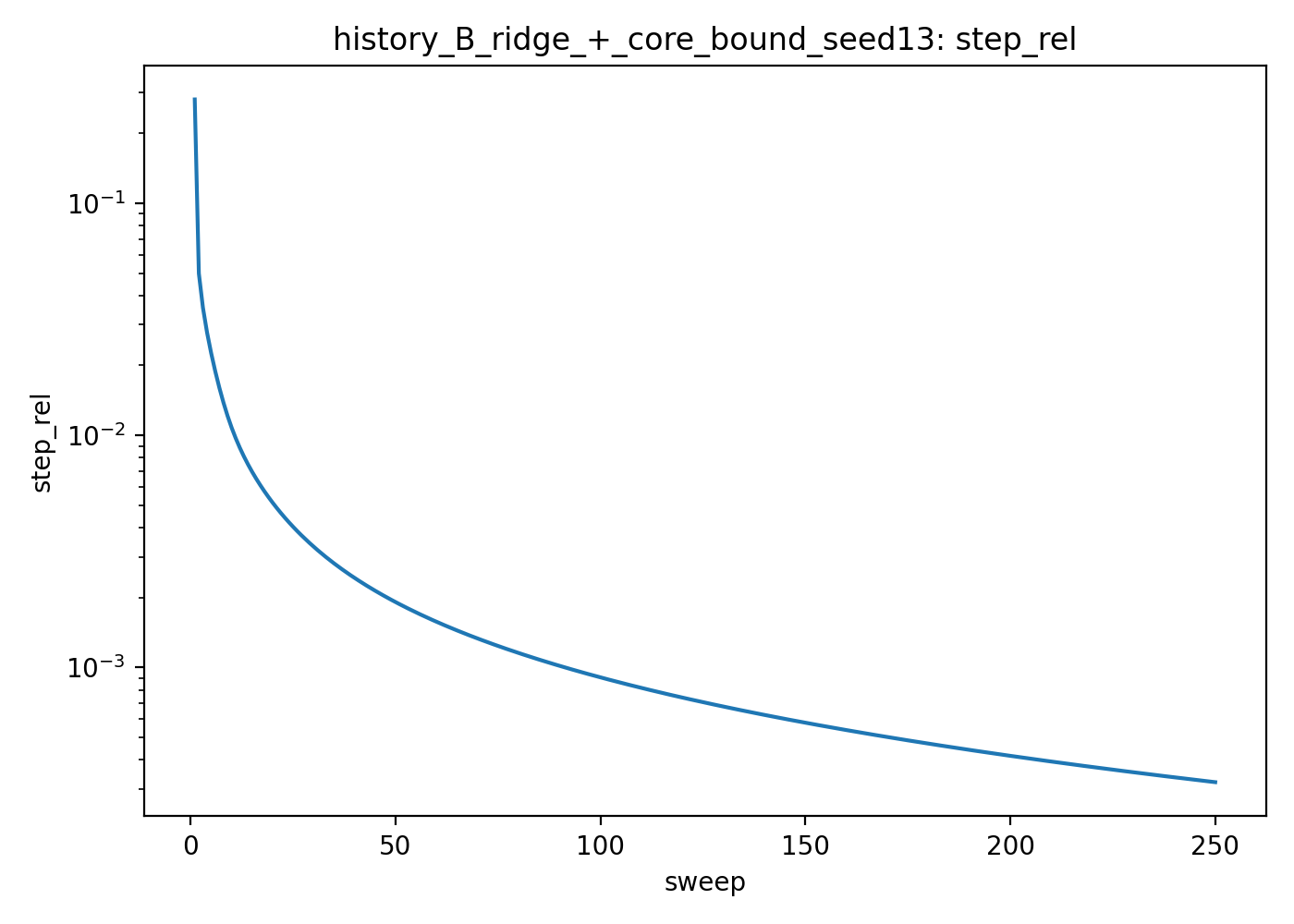}
\caption{B: relative step.}
\label{fig:main_B_step_seed13}
\end{subfigure}

\caption{Representative convergence diagnostics for ridge-only and ridge-plus-core-bound variants using the first prespecified seed. Variant A is ridge-only, while Variant B adds the core-bound mechanism. The objective and relative-step panels are the primary convergence diagnostics; residual curves are interpreted together with the regularized objective.}
\label{fig:main_convergence_A_B_seed13}
\end{figure}

Figure~\ref{fig:main_convergence_A_B_seed13} reports representative convergence histories for the ridge-only and ridge-plus-core-bound variants using the first prespecified seed. The objective and relative-step panels are the primary convergence diagnostics. The residual curves are secondary because the optimized criterion includes regularization, and a regularized update can trade a small residual change for a better-conditioned or lower-penalty representation. The histories are consistent with the descent behavior proved for the exact full-data quadratic scheme in Theorem~\ref{thm:alternating_convergence}, while the projected core-bounded implementation is treated as an algorithmic extension.

\begin{figure}[!htbp]
\centering
\begin{subfigure}[t]{0.24\textwidth}
\includegraphics[width=\linewidth]{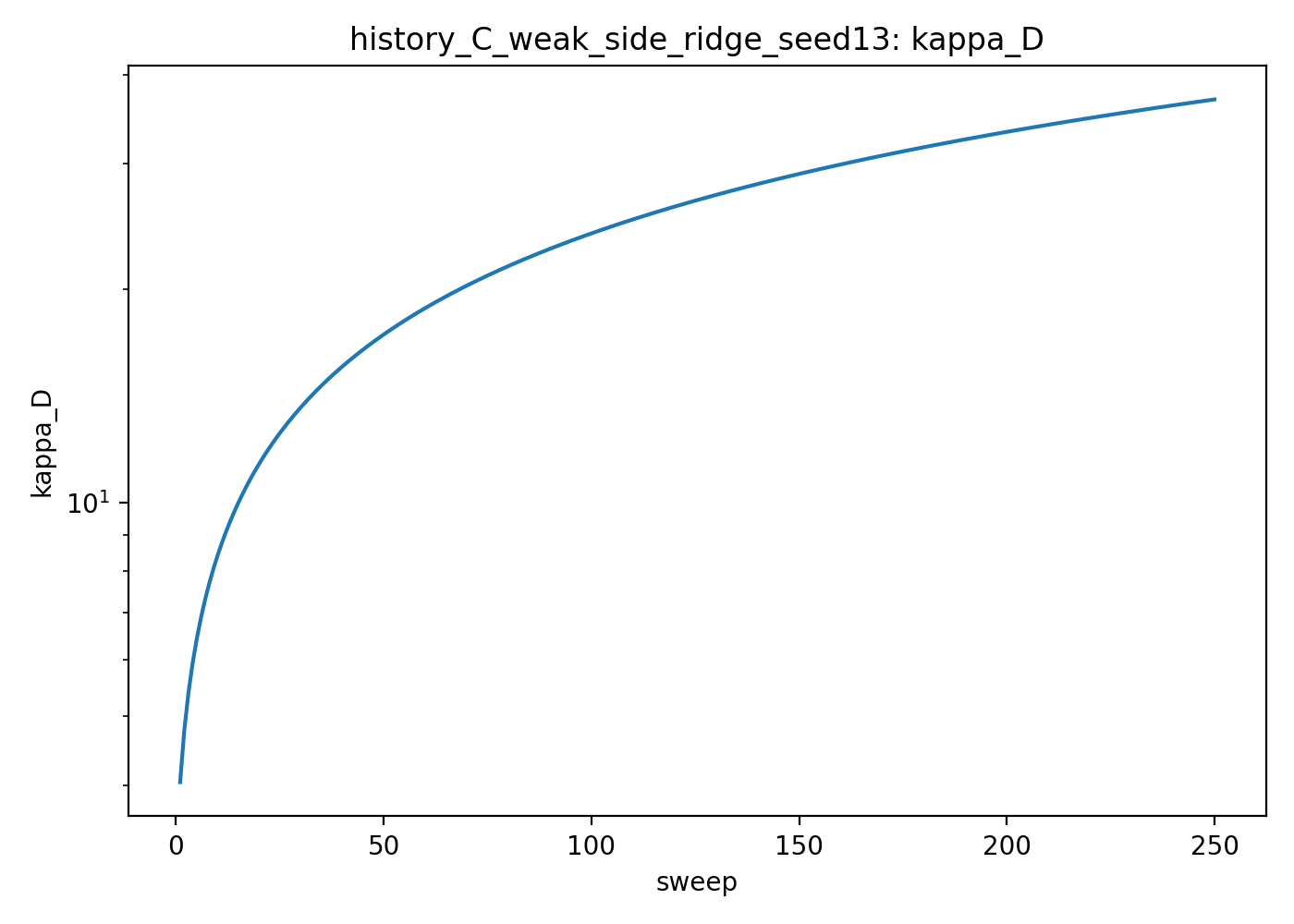}
\caption{\(\kappa(D)\).}
\label{fig:main_C_kappa_seed13}
\end{subfigure}\hfill
\begin{subfigure}[t]{0.24\textwidth}
\includegraphics[width=\linewidth]{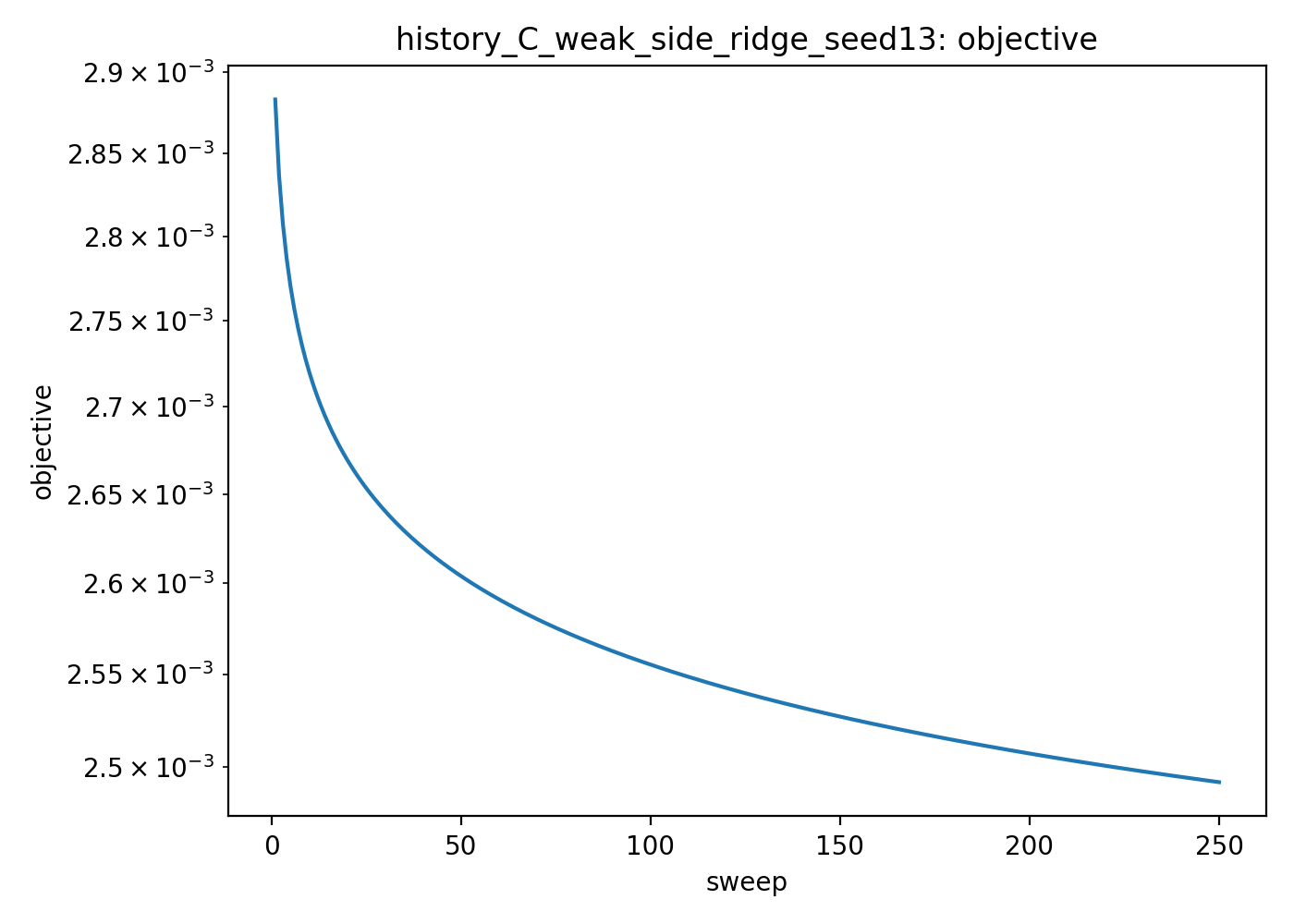}
\caption{Objective.}
\label{fig:main_C_objective_seed13}
\end{subfigure}\hfill
\begin{subfigure}[t]{0.24\textwidth}
\includegraphics[width=\linewidth]{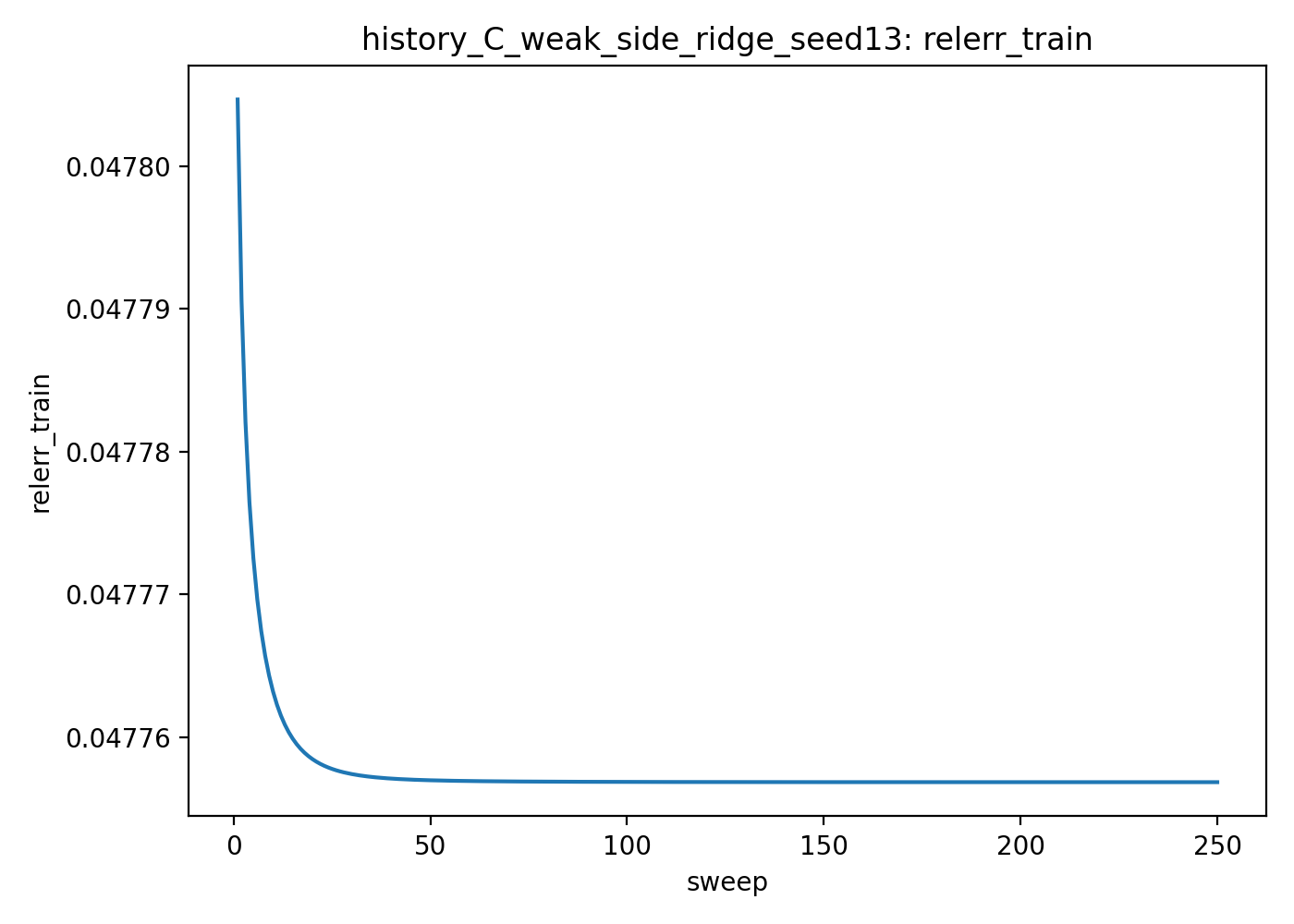}
\caption{Training residual.}
\label{fig:main_C_residual_seed13}
\end{subfigure}\hfill
\begin{subfigure}[t]{0.24\textwidth}
\includegraphics[width=\linewidth]{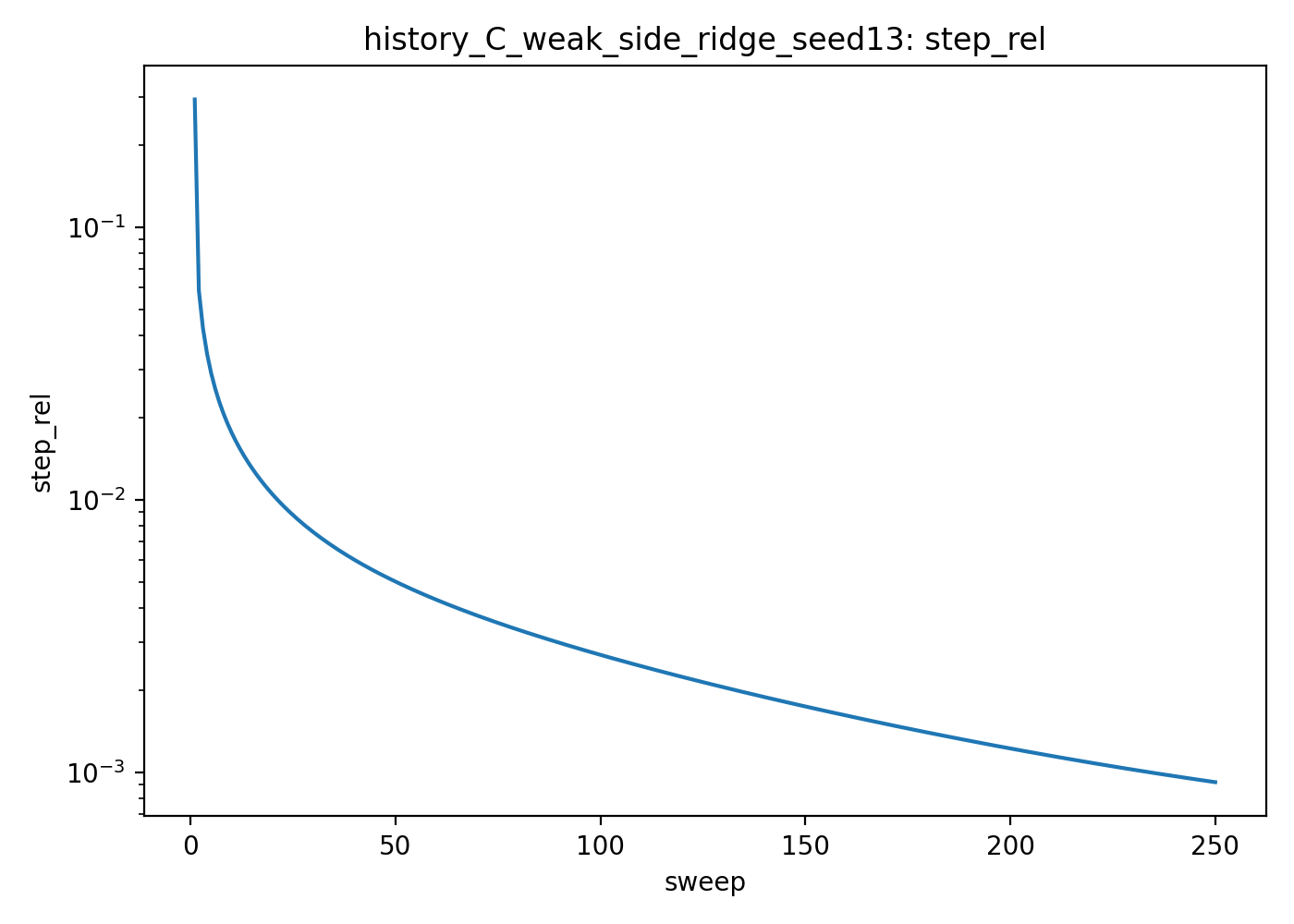}
\caption{Relative step.}
\label{fig:main_C_step_seed13}
\end{subfigure}

\caption{Representative convergence diagnostics for the weak-side-ridge ablation. This variant is included in the main text because it tests whether side-factor regularization is only cosmetic or whether it affects conditioning and convergence behavior.}
\label{fig:main_convergence_C_seed13}
\end{figure}

Figure~\ref{fig:main_convergence_C_seed13} reports the corresponding diagnostics for the weak-side-ridge ablation. This variant tests whether side-factor regularization is only cosmetic. The histories, together with Table~\ref{tab:full-ablation}, show that weakening side-factor regularization changes the conditioning behavior of the learned core.

Table~\ref{tab:full-ablation} gives the quantitative ablation results. Variant A is the ridge-only model, Variant B adds the core-bound mechanism, Variant C weakens side-factor regularization, and Variant D imposes a diagonal core. The weak-side-ridge variant achieves similar relative Frobenius error but has a much larger final \(\kappa(D)\), which supports the role of side-factor regularization in controlling scale imbalance among \(P\), \(D\), and \(Q\). The diagonal-core variant gives similar accuracy but is much slower in this implementation, so it should be viewed as a structured-core option rather than a default computational choice.
\begin{table}[!htbp]
\centering
\scriptsize
\setlength{\tabcolsep}{2.5pt}
\renewcommand{\arraystretch}{0.92}
\caption{Ablation results over the full multi-seed validation run. Variant A is ridge-only, Variant B adds the core-bound mechanism, Variant C weakens side-factor regularization, and Variant D imposes a diagonal core.}
\label{tab:full-ablation}
\resizebox{\textwidth}{!}{%
\begin{tabular}{@{}llllllll@{}}
\toprule
Variant & RelFrob & CleanErr & $\kappa(D)$ & Time (s) & Sweeps & $\rho$ & $p$ \\
\midrule
A ridge only & \(\num{0.0477306798759}\pm\num{3.01195876346e-05}\) & \(\num{0.0149505611571}\pm\num{8.46732737624e-05}\) & \(\num{5.71580368492}\pm\num{0.00854532486604}\) & \(\num{4.16795300001}\pm\num{0.125187329964}\) & \num{250} & \(\num{0.98019893136}\pm\num{3.54512362311e-05}\) & \(\num{0.996823969292}\pm\num{0.00430810912917}\) \\
B ridge + core bound & \(\num{0.0477306798759}\pm\num{3.01195876346e-05}\) & \(\num{0.0149505611571}\pm\num{8.46732737624e-05}\) & \(\num{5.71580368492}\pm\num{0.00854532486608}\) & \(\num{6.04372530004}\pm\num{0.629722347571}\) & \num{250} & \(\num{0.98019893136}\pm\num{3.54512362186e-05}\) & \(\num{0.996823969292}\pm\num{0.00430810912915}\) \\
C weak side ridge & \(\num{0.0477155534875}\pm\num{3.0107536999e-05}\) & \(\num{0.0155189521458}\pm\num{8.56485330464e-05}\) & \(\num{37.1314160986}\pm\num{0.183012551119}\) & \(\num{4.25005073999}\pm\num{0.0952295756476}\) & \num{250} & \(\num{0.981662587854}\pm\num{1.35810761008e-05}\) & \(\num{0.997893377164}\pm\num{0.00481036779912}\) \\
D diagonal core & \(\num{0.0477306794079}\pm\num{3.01197032701e-05}\) & \(\num{0.0149505977539}\pm\num{8.46582258834e-05}\) & \(\num{5.7108707724}\pm\num{0.00958419475501}\) & \(\num{201.6939074}\pm\num{2.36073268673}\) & \num{250} & \(\num{0.980198934473}\pm\num{3.54542161616e-05}\) & \(\num{0.996823200794}\pm\num{0.00430820903306}\) \\
\bottomrule
\end{tabular}%
}
\end{table}

\subsection{Runtime and storage scaling}

Table~\ref{tab:full-scaling} reports runtime and storage scaling. Randomized SVD is faster in these experiments, which is expected because it is a direct approximate spectral method. The proposed method solves repeated regularized block subproblems and records additional diagnostics, including \(\kappa(D)\), objective decay, and relative step size. The scaling results should therefore not be read as evidence that the proposed method is a faster replacement for randomized SVD. They show the trade-off: the proposed method pays additional computation to obtain a structured tri-factorization with explicit core diagnostics.
\begin{table}[!htbp]
\centering
\scriptsize
\setlength{\tabcolsep}{2.5pt}
\renewcommand{\arraystretch}{0.92}
\caption{Runtime and storage scaling over the full validation run. Randomized SVD is included as a fast spectral baseline; the proposed method reports additional core-conditioning diagnostics.}
\label{tab:full-scaling}
\resizebox{\textwidth}{!}{%
\begin{tabular}{@{}llllllll@{}}
\toprule
Dataset & Method & RelFrob & CleanErr & $\kappa(D)$ & Time (s) & Sweeps & Storage \\
\midrule
scale\_m1200\_n900 & Randomized SVD & \num{0.0196068693505} & \num{0.0039286549684} & -- & \num{0.0660234999377} & \num{0} & \num{0.0392592592593} \\
scale\_m1200\_n900 & Proposed conditioned core & \num{0.0196263432482} & \num{0.00401544734455} & \num{2.66373581197} & \num{10.0214984} & \num{120} & \num{0.0392592592593} \\
scale\_m1800\_n1200 & Randomized SVD & \num{0.0197172924924} & \num{0.00332781514597} & -- & \num{0.0493487999775} & \num{0} & \num{0.027962962963} \\
scale\_m1800\_n1200 & Proposed conditioned core & \num{0.0197366569448} & \num{0.0034301687306} & \num{2.66272997685} & \num{12.9109136} & \num{120} & \num{0.027962962963} \\
scale\_m400\_n300 & Randomized SVD & \num{0.0188257490385} & \num{0.0067447090293} & -- & \num{0.00764450000133} & \num{0} & \num{0.12} \\
scale\_m400\_n300 & Proposed conditioned core & \num{0.0188460202129} & \num{0.00677542507921} & \num{2.66656272316} & \num{2.24805250007} & \num{120} & \num{0.12} \\
scale\_m800\_n600 & Randomized SVD & \num{0.0194093293613} & \num{0.00480778406616} & -- & \num{0.0281052999198} & \num{0} & \num{0.0591666666667} \\
scale\_m800\_n600 & Proposed conditioned core & \num{0.0194289991863} & \num{0.00487533843428} & \num{2.6634401347} & \num{4.093521} & \num{120} & \num{0.0591666666667} \\
\bottomrule
\end{tabular}%
}
\end{table}

The numerical evidence supports a careful interpretation of the method. Clean low-rank recovery confirms that the proposed formulation behaves consistently with rank-constrained approximation, but it should not be presented as an improvement over SVD in the unregularized full-data problem. The noisy low-rank and ill-conditioned experiments provide the strongest evidence for the proposed framework: the method is competitive with spectral baselines while exposing core-conditioning diagnostics that those baselines do not report. The ablation results show that side-factor regularization matters for controlling the learned core. The missing-entry experiments identify a limitation of the present weighted implementation, and the scaling experiments show that randomized SVD remains faster for pure spectral compression.

Overall, the empirical evidence supports the proposed method as a regularized and diagnostically transparent tri-factorization for settings where approximation accuracy and numerical conditioning must be considered together. It does not support a claim of universal superiority over SVD, randomized SVD, or matrix-completion methods.

\subsection{Synthesis of numerical results}
\label{subsec:numerical_synthesis}

The full validation run supports a precise interpretation of the method. In the clean low-rank experiments, the spectral baselines recover the rank-constrained target to numerical precision when the target rank is sufficient, while the proposed regularized variants introduce a small but expected bias. This confirms that the method should not be presented as a replacement for truncated SVD in the unregularized full-data Frobenius problem.

The noisy low-rank experiments show that the proposed variants track the spectral baselines closely in reconstruction error and clean-target error, while also reporting a learned core condition number. This is the main empirical distinction in this regime: the proposed method gives comparable approximation together with diagnostics on the internal numerical structure of the factorization.

The ill-conditioned experiments provide the strongest evidence for the core-conditioning idea. As the spectrum becomes more severe, the unconstrained ridge core can become singular, while the conditioned-core variant enforces a finite admissible condition number. The accuracy differences are modest, but the diagnostic and stability distinction is important. This supports the central claim that the conditioning of the latent core should be treated as a controlled modeling quantity rather than an automatic consequence of the factorization.

The missing-entry experiments give a more cautious conclusion. The proposed weighted tri-factorization fits observed entries closely and becomes competitive at the highest observation rate, but it does not uniformly outperform soft-impute or ridge two-factor factorization on held-out RMSE. The weighted formulation is therefore best treated as an extension of the framework, not as the main empirical strength of the paper.

The ablation results show that side-factor regularization is not merely cosmetic. Weakening the side-factor regularization substantially increases the learned core condition number, even when the reconstruction error remains similar. The diagonal-core variant gives comparable accuracy but is much slower in the current implementation. The scaling experiments also show that randomized SVD remains faster for pure spectral compression. The proposed method therefore trades additional computation for a regularized tri-factorization with an explicit core and stability diagnostics.

\section{Discussion}
\label{sec:discussion}

The results support a focused interpretation of the proposed factorization. The representation \(A\approx PDQ\) is not a new algebraic rank class. Every product \(PDQ\) has rank at most \(r\), and every rank-\(r\) matrix admits such a representation. The useful contribution is the regularized and core-conditioned use of this parametrization. By treating \(D\) as a central core whose conditioning can be monitored, regularized, or constrained, the method separates approximation quality from the numerical quality of the compressed representation.

The theoretical results are aligned with this interpretation. Existence follows from coercive regularization, not from compactness of the unconstrained factor space. Factor uniqueness is not claimed, because invertible transformations in the latent space can leave the product \(PDQ\) unchanged while changing the individual factors. Perturbation stability is therefore stated at the product level, with factor-level stability requiring additional gauge fixing and local second-order assumptions. The convergence theorem is similarly restricted to the exact full-data quadratic alternating scheme, where the block updates are well posed and the descent argument is clean.

The numerical evidence is consistent with the intended role of the theory. The clean low-rank experiments confirm the expected SVD benchmark behavior. The noisy and ill-conditioned experiments show that the proposed variants are competitive with spectral baselines while providing additional diagnostics, especially \(\kappa(D)\), objective decay, relative step size, and the effect of regularization on the learned core. These diagnostics are not a secondary detail. They are the main reason for using a core-conditioned tri-factorization rather than only a spectral approximation.

The experiments also clarify where the method should not be overpositioned. It is not faster than randomized SVD in the scaling experiments, and it is not yet a superior matrix-completion method in the missing-entry experiments. The current evidence supports the method as a regularized, diagnostically transparent low-rank factorization for settings where approximation error and numerical conditioning must be considered together. It does not support a claim of universal superiority over classical decompositions, randomized spectral methods, or matrix-completion algorithms.

\subsection{Limitations}
\label{sec:limitations}

The numerical validation is based on controlled synthetic and diagnostic matrix regimes. This is appropriate for testing rank recovery, noise sensitivity, ill-conditioning, missing observations, convergence behavior, and ablation structure under known ground truth. It does not replace broad real-data benchmarking. Larger sparse scientific matrices, application-derived graph or kernel matrices, and real high-dimensional data matrices would be needed before making stronger claims about large-scale empirical performance.

The missing-entry formulation remains the least mature part of the current framework. In the full validation run, soft-impute and ridge two-factor factorization are stronger in several held-out RMSE comparisons. The weighted tri-factorization therefore should not be positioned as a matrix-completion competitor at this stage. Improving this case will require stronger validation-based regularization, better early stopping, and possibly a revised weighted block-update strategy.

The convergence analysis applies to the exact full-data quadratic alternating scheme. The projected core-bounded variant and the weighted missing-entry implementation are included because they are natural algorithmic extensions, but their complete convergence analysis is not established here. Extending the proof to these cases would require additional assumptions on the projection, the acceptance rule, and the weighted subproblem solvers.

The factorization is non-identifiable without additional gauge constraints. The product \(PDQ\) is the primary approximation object, while the factors \(P\), \(D\), and \(Q\) may change under invertible latent-space transformations. Factor-level interpretation therefore requires additional structure, such as normalization, diagonal or sparse cores, orthogonality, nonnegativity, or other application-specific constraints.

The current implementation also favors clarity and reproducibility over maximum speed. Randomized SVD remains faster for pure spectral compression. The proposed method is most defensible when the additional computation is justified by the need for factor regularization, an explicit core, and conditioning diagnostics.

\section{Conclusion}
\label{sec:conclusion}

This paper introduced a core-conditioned regularized matrix tri-factorization for structured low-rank approximation. The method represents a matrix by \(A\approx PDQ\), where \(D\) is a central core whose conditioning can be monitored, regularized, or constrained. The formulation is not intended to replace LU, QR, Cholesky, truncated SVD, randomized SVD, or matrix-completion methods in the settings where those methods are already the appropriate tools. Its purpose is to provide a regularized approximation framework in which reconstruction error, factor scale, and numerical conditioning can be studied together.

The theoretical analysis established the rank-constrained scope of the representation, existence under coercive regularization, non-identifiability under latent-space transformations, well-posed block updates for the quadratic full-data case, core-conditioning bounds, product-level perturbation stability, and convergence to a critical point for exact alternating minimization under standard nonconvex descent assumptions.

The full multi-seed validation supports the diagnostic value of the method. The proposed variants behave as expected on clean low-rank data, remain competitive in noisy low-rank approximation, and expose useful conditioning behavior in ill-conditioned regimes. The ablation results show that side-factor regularization affects the learned core condition number, while the scaling results show the expected computational trade-off relative to randomized SVD. The missing-entry experiments identify a limitation of the current weighted formulation and point to a clear direction for improvement.

The main contribution is therefore a careful and reusable framework for regularized matrix tri-factorization with explicit core-conditioning diagnostics. Future work should extend the convergence analysis to projected and weighted variants, improve the missing-entry formulation, test larger sparse and real-data matrices, and study structured cores such as diagonal, sparse, block-diagonal, nonnegative, or positive semidefinite cores.

\bibliographystyle{unsrtnat}
\bibliography{refs}

\appendix
\section{Supplementary Numerical Evidence and Reproducibility}
\label{app:numerical_evidence}

\paragraph{Reproducibility statement.}
All numerical results reported in the main text were generated from the full multi-seed validation pipeline accompanying this paper. The code, plotting scripts, configuration files, and raw numerical outputs are available in the public GitHub repository
\url{https://github.com/karjxenval/Core-Conditioned-Regularized-Matrix-Tri-Factorization-for-High-Dimensional-Structured-Systems}.
The repository contains the scripts used to generate the synthetic matrix families, run the clean low-rank, noisy low-rank, ill-conditioned, missing-entry, ablation, and scaling experiments, compute the convergence diagnostics, and reproduce the figures and tables reported in the paper. The raw CSV outputs are included for auditability, while the manuscript reports aggregated summary tables in the main text.

The full validation run uses the prespecified random seeds
\[
13,\ 29,\ 47,\ 71,\ 101.
\]
For each experiment, the main text reports mean values and standard deviations across repeated runs where applicable. The appendix provides additional convergence diagnostics and per-seed convergence histories that are not included in the main text in order to keep the main paper readable.

\paragraph{Scope of the supplementary evidence.}
The supplementary material is intended to support reproducibility, not to introduce additional claims beyond those made in the main text. The convergence quantities reported below are empirical diagnostics computed from the optimization histories. They should not be interpreted as universal theoretical rates. The theoretical convergence statement in Theorem~\ref{thm:alternating_convergence} applies to the exact full-data quadratic alternating scheme and guarantees convergence to a critical point under the stated assumptions.

\subsection{Aggregated empirical convergence diagnostics}
\label{app:convergence_diagnostics}

Table~\ref{tab:app-convergence-summary} reports the empirical convergence diagnostics aggregated over the per-seed histories. The quantity \(\rho\) summarizes tail objective contraction, \(q\) is the fitted sublinear decay exponent, \(p\) is the empirical fixed-point step-order estimate, \(\Delta_t\) is the final relative objective decrease, and \(s_t\) is the final relative step size. These quantities are used to document the numerical behavior of the implementation across the validation regimes.

\begin{sidewaystable}[!htbp]
\centering
\scriptsize
\setlength{\tabcolsep}{3pt}
\renewcommand{\arraystretch}{1.08}

\caption{Aggregated empirical convergence diagnostics computed from the per-seed histories. These are diagnostic quantities, not universal theoretical convergence rates.}
\label{tab:app-convergence-summary}

\resizebox{\textheight}{!}{%
	\begin{tabular}{@{}lllllll@{}}
	\toprule
	Experiment & Method & $\rho$ & $q$ & $p$ & Final $\Delta_t$ & Final $s_t$ \\
	\midrule
	ablation & A ridge only & \(\num{0.98019893136}\pm\num{3.54512362311e-05}\) & \(\num{4.80910452427}\pm\num{0.00535613769791}\) & \(\num{0.996823969292}\pm\num{0.00430810912917}\) & \(\num{7.49798617142e-08}\pm\num{7.71807418356e-11}\) & \(\num{0.000334327867522}\pm\num{2.22334308264e-05}\) \\
	ablation & B ridge + core bound & \(\num{0.98019893136}\pm\num{3.54512362186e-05}\) & \(\num{4.80910452427}\pm\num{0.0053561376982}\) & \(\num{0.996823969292}\pm\num{0.00430810912915}\) & \(\num{7.49798617138e-08}\pm\num{7.71807414288e-11}\) & \(\num{0.000334327867522}\pm\num{2.22334308266e-05}\) \\
	ablation & C weak side ridge & \(\num{0.981662587854}\pm\num{1.35810761008e-05}\) & \(\num{4.53352067054}\pm\num{0.00213828304131}\) & \(\num{0.997893377164}\pm\num{0.00481036779912}\) & \(\num{2.73853457303e-07}\pm\num{8.01115750061e-11}\) & \(\num{0.000972076663902}\pm\num{7.94409665269e-05}\) \\
	ablation & D diagonal core & \(\num{0.980198934473}\pm\num{3.54542161616e-05}\) & \(\num{4.80910393665}\pm\num{0.00535664772219}\) & \(\num{0.996823200794}\pm\num{0.00430820903306}\) & \(\num{7.49799933277e-08}\pm\num{7.71481493802e-11}\) & \(\num{0.000334347562159}\pm\num{2.22156656592e-05}\) \\
	clean\_low\_rank & Proposed ridge & \(\num{0.978913129026}\pm\num{2.13330262158e-08}\) & \(\num{5.05067398279}\pm\num{4.54884587242e-06}\) & \(\num{0.9969357381}\pm\num{5.56573169058e-08}\) & \(\num{1.39798941103e-07}\pm\num{1.00454538031e-13}\) & \(\num{0.00046364562816}\pm\num{7.40892501026e-10}\) \\
	clean\_low\_rank & Proposed conditioned core & \(\num{0.978913129026}\pm\num{2.1333068974e-08}\) & \(\num{5.05067398279}\pm\num{4.54884567542e-06}\) & \(\num{0.9969357381}\pm\num{5.56572941159e-08}\) & \(\num{1.39798941103e-07}\pm\num{1.00454491789e-13}\) & \(\num{0.00046364562816}\pm\num{7.40892560714e-10}\) \\
	ill\_conditioned & Proposed ridge & \(\num{0.980474604856}\pm\num{0.000518400895401}\) & \(\num{4.75619457389}\pm\num{0.100920290488}\) & \(\num{0.99559715531}\pm\num{0.00124515114083}\) & \(\num{1.12791477314e-07}\pm\num{5.39908123572e-08}\) & \(\num{0.000315129576077}\pm\num{2.82972742874e-05}\) \\
	ill\_conditioned & Proposed conditioned core & \(\num{0.980474604856}\pm\num{0.000518400895398}\) & \(\num{4.75619457389}\pm\num{0.100920290489}\) & \(\num{0.840278603415}\pm\num{0.196918367525}\) & \(\num{1.12791477315e-07}\pm\num{5.39908123591e-08}\) & \(\num{0.000449669236921}\pm\num{0.000191131135443}\) \\
	missing\_entries & Soft-impute & \(\num{0.975883867714}\pm\num{0.0039872931356}\) & \(\num{5.44599452433}\pm\num{0.702140481862}\) & \(\num{0.991139847311}\pm\num{0.00520827887292}\) & \(\num{1.18251228098e-06}\pm\num{1.45810070077e-06}\) & \(\num{0.000675695946018}\pm\num{0.000494490116798}\) \\
	missing\_entries & Ridge two-factor & \(\num{0.978382163144}\pm\num{0.00211517419266}\) & \(\num{5.19260640779}\pm\num{0.373610242355}\) & \(\num{0.99518769401}\pm\num{0.00443742356245}\) & \(\num{4.39307313659e-07}\pm\num{3.68556592479e-07}\) & \(\num{0.00215662544998}\pm\num{0.00106415046054}\) \\
	missing\_entries & Proposed ridge & \(\num{0.95844231195}\pm\num{0.00535563663987}\) & \(\num{4.92151354469}\pm\num{0.366712425282}\) & \(\num{0.967903126046}\pm\num{0.0416393165856}\) & \(\num{2.52843056387e-06}\pm\num{1.50667685871e-06}\) & \(\num{0.00420546989908}\pm\num{0.00174014200685}\) \\
	missing\_entries & Proposed conditioned core & \(\num{0.958442311903}\pm\num{0.00535563671053}\) & \(\num{4.92151354189}\pm\num{0.366712421167}\) & \(\num{0.967882348853}\pm\num{0.0416767351118}\) & \(\num{2.52843056387e-06}\pm\num{1.50667685871e-06}\) & \(\num{0.00420602711284}\pm\num{0.00174028497562}\) \\
	noisy\_low\_rank & Ridge two-factor & \(\num{0.981997219715}\pm\num{0.000816926903628}\) & \(\num{4.46560410755}\pm\num{0.154115994972}\) & \(\num{0.997069144145}\pm\num{0.00150533944814}\) & \(\num{1.5457563109e-06}\pm\num{8.72888929677e-07}\) & \(\num{0.000591391978422}\pm\num{0.000245181503165}\) \\
	noisy\_low\_rank & Proposed ridge & \(\num{0.982071061258}\pm\num{4.20288176981e-05}\) & \(\num{4.44880261472}\pm\num{0.00788043836111}\) & \(\num{0.99869557108}\pm\num{5.12712331443e-05}\) & \(\num{1.74609766761e-07}\pm\num{1.89589558479e-09}\) & \(\num{0.000206067903525}\pm\num{3.67228076951e-06}\) \\
	noisy\_low\_rank & Proposed conditioned core & \(\num{0.982071061257}\pm\num{4.20288169134e-05}\) & \(\num{4.44880261472}\pm\num{0.007880438358}\) & \(\num{0.99869557108}\pm\num{5.12712331996e-05}\) & \(\num{1.74609766758e-07}\pm\num{1.89589558437e-09}\) & \(\num{0.000206067903525}\pm\num{3.67228076954e-06}\) \\
	scaling & Proposed conditioned core & \(\num{0.959382677095}\pm\num{6.27022610883e-07}\) & \(\num{4.68489771813}\pm\num{5.59614068924e-05}\) & \(\num{0.990676046207}\pm\num{5.16295208523e-06}\) & \(\num{5.78844704231e-07}\pm\num{1.92767662278e-10}\) & \(\num{0.00116621696642}\pm\num{4.27213530156e-08}\) \\
	\bottomrule
	\end{tabular}
	}

\end{sidewaystable}

\subsection{Additional per-seed convergence histories}
\label{app:additional_convergence_figures}

The main text reports representative convergence histories using seed \(13\), the first prespecified seed. This subsection reports the remaining per-seed histories. These figures are included to show that the displayed convergence behavior is not tied to a single representative run. Each history records the learned core condition number \(\kappa(D)\), the objective value, the relative training residual, and the relative step size across alternating sweeps.

The figures are grouped by ablation variant. Variant A denotes the ridge-only formulation, Variant B denotes ridge regularization with the core-bound mechanism, Variant C weakens side-factor regularization, and Variant D imposes a diagonal core. The same diagnostics are reported across variants to make the convergence behavior and conditioning effects directly comparable.

\end{document}